\documentclass[leqno, 12pt]{article}

\usepackage[latin1]{inputenc}
\usepackage{amsmath,amssymb}
\usepackage{latexsym}

\newtheorem{theorem}{Theorem}[section]

\newtheorem{lemma}[theorem]{Lemma}

\numberwithin{equation}{section}

\def \be{\begin{equation}}
\def \ee{\end{equation}}
\def \bt{\begin{theorem}}
\def \et{\end{theorem}}
\def \bea{\begin{eqnarray}}
\def \eea{\end{eqnarray}}
\def \bas{\begin{eqnarray*}}
\def \eas{\end{eqnarray*}}

\def \al{\beta}
\def \bb{\beta}
\def \ga{\gamma}
\def \Ga{\Gamma}
\def \de{\delta}

\def \ep{\epsilon}

\def \la{\lambda}
\def \La{\Lambda}

\def \si{\sigma}
\def \Si{\Sigma}

\def \th{\theta}

\def \ff{\infty}
\def \wh{\widehat}
\def \wt{\widetilde}

\def \rar{\rightarrow}

\def \R{{\bf R}}

\def \Z{{\bf Z}}

\def \({\left(}
\def \){\right)}

\def \lc{\left\{}
\def \rc{\right\}}

\def \nn{\nonumber}
\def \Proof{{\bf Proof: }}

\def \bc{\begin{center} }
\def \ec{\end{center} }
\def \bs{\begin{slide} }
\def \es{\end{slide} }

\def\R{{\mathbb R}}
\def\Z{{\mathbb Z}}

\def\lam{\lambda}
\def\al{\beta}

\def\phi{\varphi}

\def\proof{{\medskip\noindent {\bf Proof. }}}

\def\qed{{\hfill $\square$ \bigskip}}

\def\square{{\vcenter{\vbox{\hrule height.3pt
          \hbox{\vrule width.3pt height5pt \kern5pt
             \vrule width.3pt}
          \hrule height.3pt}}}}

 \def\bE {{\Bbb E}}

\def\bP {{\Bbb P}}

\def\square{{\vcenter{\vbox{\hrule height.3pt
           \hbox{\vrule width.3pt height5pt \kern5pt
              \vrule width.3pt}
           \hrule height.3pt}}}}
\def\qed{{\hfill $\square$ \bigskip}}
\def\ol{\overline}
\def\E{{\bE}}
\def\P{{\bP}}

\begin{document}

\title{Large Deviations for  Riesz Potentials of Additive Processes}
\author{Richard Bass\thanks{Research partially supported by NSF grant
\#DMS-0601783}
     \, \, Xia Chen\thanks {Research   partially supported by NSF grant
\#DMS-0405188.}\, \, Jay Rosen\thanks {Research partially supported  by grants
from the NSF  and from PSC-CUNY.}}

\maketitle

\bibliographystyle{amsplain}

\begin{abstract} We study functionals of the form
\[\zeta_{t}=\int_0^{t}\cdots\int_0^{t}
\vert X_1(s_1)+\cdots+ X_p(s_p)\vert^{-\sigma}ds_1\cdots ds_p\] 
\end{abstract}
where $X_1(t),\cdots, X_p(t)$ are i.i.d. $d$-dimensional  symmetric stable processes of
index $0<\bb\le 2$. We obtain results about the large deviations and laws of the iterated logarithm
 for $\zeta_{t}$.

\section{Introduction }

Let $X_1(t),\cdots, X_p(t)$ be i.i.d. $d$-dimensional symmetric stable process of
index $0<\bb\le 2$. We use the notation $X(t)$ for a stable process with the same
distribution as $X_1(t),\cdots X_p(t)$. Thus
\begin{equation}
\E e^{ i\lambda \cdot X_{ t}}=e^{ -t| \lambda|^{ \bb}}\hskip.2in t\ge 0,
\hskip.1in
\lambda\in\R^d.\label{a1.0}
\end{equation} 
In this paper we study
\begin{equation}
\zeta([0,t_1]\times\cdots\times [0,t_p])=\int_0^{t_1}\cdots\int_0^{t_p}
\vert X_1(s_1)+\cdots+ X_p(s_p)\vert^{-\sigma}ds_1\cdots ds_p\label{a1.1}
\end{equation} 
and more generally 
\begin{equation}
\zeta^{z}([0,t_1]\times\cdots\times [0,t_p])=\int_0^{t_1}\cdots\int_0^{t_p}
\vert X_1(s_1)+\cdots+ X_p(s_p)-z\vert^{-\sigma}ds_1\cdots ds_p\label{a1.1}
\end{equation}
for $z\in R^{d}$.
We show below that $\zeta^{z}([0,t_1]\times\cdots\times [0,t_p])$ is
finite almost surely if
\begin{equation} 0<\sigma <\min\{p\bb, d\}.\label{a1.2}
\end{equation}
The random field 
$\bar X(t_{1},\ldots,t_{p})=X_1(t_{1})+\cdots + X_p(t_{p})$ is known as an additive process, and its occupation measure $\mu_{A}$  for  $A\in R_{+}^{p}$ is the measure on $R^{d}$  defined by
\begin{equation}
\mu_{A}(B)=\int_{A}1_{\{X_1(s_{1})+\cdots + X_p(s_{p})\in B\}}\,ds_{1}\cdots\,ds_{p}. \label{a1.89}
\end{equation}
With this notation we have
\begin{equation}
\zeta^{z}([0,t_1]\times\cdots\times [0,t_p])=\int {1 \over |x-z|^{-\si}}\,\,\mu_{[0,t_1]\times\cdots\times [0,t_p]}(\,dx)\label{a1.88}
\end{equation}
so that $\zeta^{z}([0,t_1]\times\cdots\times [0,t_p])$ is the Riesz potential of the occupation measure $\mu_{[0,t_1]\times\cdots\times [0,t_p]}$. (In the terminology of \cite{D}, $\zeta^{z}([0,t_1]\times\cdots\times [0,t_p])$ is the Riesz-Frostman potential of the occupation measure.)

Because they locally resemble stable sheets, but are 
more amenable to analysis, additive stable processes first 
arose to simplify the study of stable sheets (see
Dalang and Walsh \cite{DalangWalsha,DalangWalshb}, Kahane \cite{Kahane} and Kendall \cite{Kendall}).
They also arise in the theory of intersections and self 
intersections of stable processes (see Le Gall, Rosen and
Shieh \cite{LGR}, Fitzsimmons and Salisbury \cite{FitzsimmonsSalisbury}, Khoshnevisan and Xiao \cite{KhoshnevisanXiao}).
In addition, the study of additive processes has connections with 
probabilistic potential theory. We  refer the reader to  Hirsch and Song \cite{HirschSong}, Khoshnevisan
\cite{Khoshnevisan}, Khoshnevisan and Shi \cite{KhoshnevisanShi}, Khoshnevisan and Xiao \cite{KhoshnevisanXiao}
for detailed discussion and further references. The present paper is a direct outgrowth of \cite{Cadditive}.

We are interested in Riesz potentials for two reasons. 
 First of all, they provide an opportunity to study functionals of the paths which  are, in some ways, more singular than local times. More precisely, although Riesz potentials involve the functions $|x|^{-\si}$ while local times involve the more singular delta `function', much of our analysis in both cases involves Fourier transforms, and the Fourier transform of $\de_{0}$  is $1$, while that of $|x|^{-\si}$ is $c|x|^{-(d-\si)}$. 
The second reason involves generalizations of the polaron problem.  
Donsker and Varadhan 
\cite{DV} show that for Brownian motion in $R^{3}$
\begin{eqnarray} &&
\lim_{t\to\infty}{1\over t}\log E\exp\bigg\{{1 \over t}  
\int_0^{ t}\int_0^{ t} {1\over |X_s-X_r|
}\,dr\,ds \bigg\}\label{dv}\\ &&
= \sup_{g\in \mathcal{F}_2}
\bigg\{\int_{R^3}\int_{R^3}
{ g^2(x) g^2 (y)\over |x-y|}dx\,dy -
{ 1\over 2}\| \nabla f \|_{ 2}^{ 2}\bigg\}.\nonumber
\end{eqnarray}
The object in the exponential involves a Riesz potential but here we have a single process as opposed to several independent processes.

\bt\label{theo-jc} Under (\ref{a1.2}),  
$\zeta^{z}([0,t_1]\times\cdots\times [0,t_p])$ is jointly continuous in 
$z,t_{1},\ldots, t_{p}$,  almost surely.
\et

We note for later reference that by scaling we have
\begin{equation}
\zeta^{z}([0,t]^p)\buildrel d\over =t^{p\bb -\sigma\over\bb}
\zeta^{z/t^{1/\bb}}([0,1]^p).\label{a1.3}
\end{equation}

For $0<\si<d$ let
\begin{equation}
\phi_{ d-\si}(\lambda)={C_{ d,\si} \over |\la|^{ d-\si}}\label{a1.4k}
\end{equation}
where $C_{ d,\si}=\pi^{-d/2}2^{ -\si}\Ga( { d-\si\over 2})/\Ga( { \si\over
2})$.
 Write
\begin{equation}
\rho=\sup_{\vert\vert f\vert\vert_2=1}
\int_{\R^d}\bigg[\int_{\R^d}{f(\lambda +\gamma)f(\gamma)\over\sqrt{
1+\psi(\lambda +\gamma)}\sqrt{ 1+\psi(\gamma)}}d\gamma\bigg]^p\phi_{
d-\si}(\lambda)d\lambda\label{a1.7}
\end{equation} where $\psi(\lambda)=| \lambda|^{ \bb}$ is the characteristic
exponent of the stable processes. Clearly, $\rho >0$.  We will prove below that $\rho <\infty$ under
condition (\ref{a1.2}). 
 
 Our main theorem is the large deviation principle for
$\zeta\big([0,t]^p\big)$. By the scaling property (\ref{a1.3}) 
 we need only consider $\zeta\big([0,1]^p\big)$
 in the following theorem.

\bt\label{theo-1} Under (\ref{a1.2}),
\begin{equation}
\lim_{t\to\infty}t^{-\bb/\sigma}\log\P\Big\{
 \zeta([0,1]^p)\ge t\Big\}=- {\sigma\over\bb}
\Big({p\bb -\sigma\over p\bb}
\Big)^{p\bb -\sigma\over\sigma}\rho^{-\bb/\sigma}\label{a1.28}
\end{equation} where $\rho$ is given in (\ref{a1.7}).
\et

The next Theorem treats the large deviations of \[\zeta^{\ast}\big([0,1]^p\big)=:\sup_{z\in R^{d}}\zeta^{z}\big([0,1]^p\big).\]

 \bt\label{theo-sup} Under (\ref{a1.2}), when $\bb=2$
\begin{equation}
\lim_{t\to\infty}t^{-\bb/\sigma}\log\P\Big\{
\zeta^{\ast}([0,1]^p)\ge t\Big\}=- {\sigma\over\bb}
\Big({p\bb -\sigma\over p\bb}
\Big)^{p\bb -\sigma\over\sigma}\rho^{-\bb/\sigma}\label{a1.28sup}
\end{equation} 
while for $\bb<2$, for some $0<C<\ff$
\begin{eqnarray}
&&- {\sigma\over\bb}
\Big({p\bb -\sigma\over p\bb}
\Big)^{p\bb -\sigma\over\sigma}\rho^{-\bb/\sigma}\leq 
\liminf_{t\to\infty}t^{-\bb/\sigma}\log\P\Big\{
\zeta^{\ast}([0,1]^p)\ge t\Big\}
\label{a1.28supb}\\
&&  \hspace{1.5 in} \leq \limsup_{t\to\infty}t^{-\bb/\sigma}\log\P\Big\{
\zeta^{\ast}([0,1]^p)\ge t\Big\}\leq -C\nonumber
\end{eqnarray}
where $\rho$ is given in (\ref{a1.7}).
\et
We believe that (\ref{a1.28sup}) holds for all $\bb$.

We can also find a law of the iterated logarithm for $\zeta^{z}([0,t]^p)$ and $\zeta^{\ast}([0,t]^p)$.
 \bt\label{theo-lil} Under (\ref{a1.2}),  
\begin{equation}
\limsup_{t\to\infty}t^{-{p\bb-\sigma \over \bb}}(\log\log t)^{-\si/\bb}
\zeta ([0,t]^p)=\({\sigma\over\bb}\)^{-\si/\bb}
\Big({p\bb -\sigma\over p\bb}
\Big)^{\sigma-p\bb \over\bb}\rho \label{a1.28lil}
\end{equation}
almost surely and when $\bb=2$
\begin{equation}
\limsup_{t\to\infty}t^{-{p\bb-\sigma \over \bb}}(\log\log t)^{-\si/\bb}
\zeta^{\ast}([0,t]^p)=\({\sigma\over\bb}\)^{-\si/\bb}
\Big({p\bb -\sigma\over p\bb}
\Big)^{\sigma-p\bb \over\bb}\rho. \label{a1.28lils}
\end{equation}
\et

We can obtain a variational expression for $\rho.$   Let $\bb\leq 2$ and set
\begin{equation}
 \mathcal{E}_{ \bb}( f,f)=:(2\pi)^{-d}\int_{ R^d} |\la|^{ \bb}|\wh{f}( \la)|^2\,d\la.
\label{7.0a}
\end{equation} Let
\begin{equation}
\mathcal{F}_{ \bb}=\{f\in L^2( R^d)\,|\,\|f\|_{ 2}=1\,,\,
\,\mathcal{E}_{ \bb}( f,f)<\ff\}.\label{7.0b}
\end{equation}
We show below that under condition (\ref{a1.2})
\begin{equation}
\La_{\si}=:\sup_{g\in{\cal F}_\beta}\bigg\{ 
\bigg(\int_{ (R^{ d})^{ p}} {\prod_{ j=1}^{ p}g^{ 2}( x_{ j})\over |x_{ 1}+\cdots+x_{
p}|^{ \si}}
\prod_{ j=1}^{ p}\,d x_{ j}\bigg)^{1/p}-{\cal E}_\bb (g,g)\bigg\}<\ff.\label{67.0v}
\end{equation}

\bt\label{theo-varia} Under condition (\ref{a1.2})
\begin{equation}
\rho=( 2\pi)^{ -d}( \La_{\si})^{p- \si/\bb}.\label{b6.9wj}
\end{equation}
\et

 We now prove that $\rho <\infty$ under 
condition (\ref{a1.2}). This will follow from the next Lemma and the fact 
 that $\bb pd/\si>d$ by (\ref{a1.2}).
 
\begin{lemma}\label{lem-sob} For any $f,g,h$ with $h\geq 0$
\begin{equation}
\(\int_{\R^d}\bigg[\int_{\R^d}{|f(\lambda +\gamma)g(\gamma)|\over\sqrt{
h(\lambda +\gamma)}\sqrt{h(\gamma)}}d\gamma\bigg]^p\phi_{
d-\si}(\lambda)d\lambda\)^{ 1/p}\leq C\|  f  \|_2 \,\|  g  \|_2\,\| h^{-1}  \|_{
pd/\si}.\label{sob.1}
\end{equation}
\end{lemma}

{\bf  Proof of Lemma \ref{lem-sob}}  By H\"older's
inequality
\begin{eqnarray} &&
\bigg[
\int_{\R^d}{|f(\lambda +\gamma)g(\gamma)|\over\sqrt{
h(\lambda +\gamma)}\sqrt{h(\gamma)}}d\gamma\bigg]^p\label{a1.21}\\ &&
=\bigg[
\int_{\R^d}|f(\lambda +\gamma)g(\gamma)|^{( p-1)/p}{|f(\lambda
+\gamma)g(\gamma)|^{ 1/p}\over\sqrt{
h(\lambda +\gamma)}\sqrt{h(\gamma)}}d\gamma\bigg]^p\nn\\ &&
 \le\bigg(\int_{\R^d}\vert f(\lambda +\gamma)g(\gamma)\vert
d\gamma\bigg)^{p-1}
\int_{\R^d}{\vert f(\lambda +\gamma)g(\gamma)\vert \over
\big(h(\lambda +\gamma)\big)^{p/2}
\big(h(\gamma)\big)^{p/2}}d\gamma.\nonumber
\end{eqnarray} By the Cauchy-Schwartz inequality and translation invariance,
\begin{equation}
\int_{\R^d}\vert f(\lambda +\gamma)g(\gamma)\vert d\gamma\le
\|f\|_{ 2}\|g\|_{ 2}.\label{a1.22}
\end{equation} Hence,
\begin{eqnarray} &&
\int_{\R^d}\bigg[\int_{\R^d}{|f(\lambda +\gamma)g(\gamma)|\over\sqrt{
h(\lambda +\gamma)}\sqrt{ h(\gamma)}}d\gamma\bigg]^p  \phi_{
d-\si}(\lambda)\,d\lambda\label{a1.23}\\ && 
\le \|f\|^{ p-1}_{ 2}\|g\|^{ p-1}_{ 2}\int_{\R^d}
\bigg(\int_{\R^d}{\vert f(\lambda +\gamma)g(\gamma)\vert \over
\big(h(\lambda +\gamma)\big)^{p/2}
\big(h(\gamma)\big)^{p/2}}d\gamma\bigg)  \phi_{
d-\si}(\lambda)\,d\lambda\nonumber
\\ &&  =C_{ d,\si} \|f\|^{ p-1}_{ 2}\|g\|^{ p-1}_{ 2} \int_{\R^d}\int_{\R^d}{F( \ga)G(
\la)\over |\la-\ga|^{ d-\si}}
 d\gamma  \,d\lambda\nonumber
\end{eqnarray} where
\begin{equation} F( \ga)=:{|f(\gamma)| \over
\big(h(\gamma)\big)^{p/2}},\,\hspace{.2in}G( \la)=:{|g(\la)| \over
\big(h(\la)\big)^{p/2}}.\label{a1.25}
\end{equation} 
Sobolev's inequality, \cite[p. 275]{D}, says that
\begin{equation}
\int_{\R^d}\int_{\R^d}{F( \ga)G( \la)\over |\la-\ga|^{ d-\si}}
 d\gamma  \,d\lambda\leq C\|F\|_{r}\|G\|_{s
}\label{sobolineq}
\end{equation}
for any $r,s>1$ with $s^{-1}+r^{-1}=1+\si/d$. 
In particular,
\begin{equation}
\int_{\R^d}\int_{\R^d}{F( \ga)G( \la)\over |\la-\ga|^{ d-\si}}
 d\gamma  \,d\lambda\leq C\|F\|_{2d/(d+\si) }\|G\|_{2d/(d+\si)
}\label{a1.26}
\end{equation} and by H\"older's inequality
\begin{eqnarray} &&
\int |F( \ga)|^{2d/(d+\si) }\,d\la=\int {|f(\gamma)|^{2d/(d+\si) } \over
\big(h(\gamma)\big)^{pd/(d+\si) }}\,d\la\label{a1.2x}\\&&
\le \|    |f|^{2d/(d+\si) }  \|_{(d+\si)/d }\,\,\|    h^{-pd/(d+\si) } 
\|_{ (d+\si)/\si}\nn\\&&
\le \|  f  \|_2^{2d/(d+\si) }\,\,\|   h^{-pd/\si}  \|_1^{
\si/(d+\si)}.\nn
\end{eqnarray} Thus 
\begin{eqnarray} &&
\|F\|_{2d/(d+\si) }\leq  
\|  f  \|_2\,\,\|   h^{-1}  \|_{ pd/\si}^{ p/2}\label{a1.2y}
\end{eqnarray} 
Our Lemma follows.\qed

We next show that $\zeta^{z}([0,t_1]\times\cdots\times [0,t_p])$ is
finite almost surely under 
condition (\ref{a1.2}). 
Let $p_{ t}( x)$ denote the transition density for the symmetric stable process in
$R^{ d}$ of index $\bb$. As usual, we define the $\al$-potential density by
\begin{equation}
u^{ \al}( x)=\int_{ 0}^{ \ff}e^{ -\al t}p_{ t}( x)\,dt.
\label{u.2}
\end{equation}
By independence  
\begin{eqnarray} &&
\E\(\zeta^{z}([0,t_1]\times\cdots\times [0,t_p])\)\label{u.3}\\ &&
=\int_0^{t_1}\cdots\int_0^{t_p}\int {1 \over |x_{1}+\cdots+x_{p}-z|^{\si}}\prod_{j=1}^{p}
p_{s_{j}}(x_{j})
\,dx_j\, ds_j\nonumber\\ &&
\leq e^{\sum_{j=1}^{p}t_{j}}\int {1 \over |x_{1}+\cdots+x_{p}-z|^{\si}}\prod_{j=1}^{p}
\int_0^{t_p}e^{-s_{j}}p_{s_{j}}(x_{j})\, ds_j\,dx_j\nn\\ &&\leq e^{\sum_{j=1}^{p}t_{j}}\int {1 \over |x_{1}+\cdots+x_{p}-z|^{\si}}
\prod_{j=1}^{p} u^{1}(x_{j})\,dx_j\nn\\ &&\leq e^{\sum_{j=1}^{p}t_{j}}\int {1 \over |x-z|^{\si}}
(u^{1}\ast\cdot\ast u^{1})(x)\,dx\nn
\end{eqnarray} 
where $(u^{1}\ast\cdot\ast u^{1})$ is the $p$-fold convolution of $u^{1}$ with itself. 
$u^{ 1}( x)$ is integrable, monotone decreasing in $|x|$,
 and asymptotic at $x=0$ to  $u^{ 0}( x)=C|x|^{-\max (0, (d-\bb))}$. Hence $(u^{1}\ast\cdot\ast u^{1})$ is integrable 
and bounded except (possibly) at $x=0$ where it is  asymptotic to  $C|x|^{-\max (0, (d-p\bb))}$. 
Hence (\ref{u.3}) is finite if (\ref{a1.2}) holds.\qed

Outline: In Section \ref{sec-kill} we prove Theorem \ref{theo-jc} and provide the general outline for our proof of the main result of this paper,  Theorem \ref{theo-1}, on large deviations. The details are carried out in Sections \ref{sec-lb}-\ref{sec-Mlim}. Section \ref{sec-varia} is devoted to the proof of the variational formula of Theorem \ref{theo-varia}, while in Section \ref{sec-sup} we prove Theorem \ref{theo-sup} on large deviations for  $\zeta^{\ast}$. Section \ref{sec-lil} establishes Theorem \ref{theo-lil} on laws of the iterated logarithm.  Finally,  an Appendix, Section \ref{sec-app},  provides certain Sobolev-type inequalities which are needed for our proofs.

Conventions: We define
\begin{equation}
\wh{f}( \la)=\int_{\R^d}e^{ ix\cdot\la}f( x)\,dx.\label{c6.1}
\end{equation}
With this notation
\begin{equation}
f( x)=( 2\pi)^{ -d}\int_{\R^d}e^{ -ix\cdot\la}\wh{f}( \la)\,dx,\label{c6.2}
\end{equation}
\begin{equation}
\wh{f\ast g}( \la)=\wh{f}( \la)\wh{g}( \la),\hspace{ .4in}\wh{f g}( \la)
=( 2\pi)^{ -d}\wh{f}( \la)\ast
\wh{g}(\la),\label{c6.3}
\end{equation}
and Parseval's identity is 
\begin{equation}
(f,g)_{ 2}=( 2\pi)^{ -d}(\wh{f}( \la),\wh{g}( \la))_{ 2}.\label{c6.4}
\end{equation}
If $\Phi\in \mathcal{S'}(R^{d})$, the set of tempered distributions on $R^{d}$, we use 
$\mathcal{F}(\Phi)$ to denote the Fourier transform of $\Phi$, so that for any 
$f\in \mathcal{S}(R^{d})$
\begin{equation}
\mathcal{F}(\Phi)(f)=\Phi (\wh{f}).\label{a1.4j}
\end{equation}
It is well known. e.g. \cite[p. 156]{D}, that $\phi_{ d-\si}\in \mathcal{S'}(R^{d})$ for any $0<\si<d$  
 and 
\begin{equation}
\mathcal{F}(\phi_{ d-\si})={1 \over |x|^{\si}}.\label{a1.4}
\end{equation}

\section{Killing at exponential times}\label{sec-kill}

We begin by citing  \cite[Lemma 2.3]{KM}.

\begin{lemma}\label{lem-1.4}  Let $Y$
be any non-negative random variable and let $\theta>0$ be fixed. Assume that
\begin{equation}
\lim_{n\to\infty}{1\over n}\log {1\over (n!)^\theta}\E Y^n=-\kappa\label{km.1}
\end{equation} for some $\kappa\in\R$. Then we have
\begin{equation}
\lim_{t\to\infty} t^{-1/\theta}\log \P\{Y\ge t\}=-\theta e^{\kappa/\theta}.
\label{km.2}
\end{equation}
\end{lemma}

In \cite{KM},  K\"onig and M\"orters assume that $\theta$ is a  positive
integer. By examining their proof, we find that $\theta$ can be any positive number.

Using this Lemma, Theorem \ref{theo-1} will follow from
\begin{equation}
\lim_{m\to\infty}{1\over m}\log {1\over (m!)^{\sigma/\bb}}\E\zeta([0,1]^p)^m
=\log\Big({p\bb \over p\bb-\sigma}\Big)^{p\bb-\sigma\over \bb} +\log
 \rho .\label{a1.8}
\end{equation}

In this section we show that (\ref{a1.8}) follows from
\begin{equation}
\lim_{m\to\infty}{1\over m}\log {1\over (m!)^p}
\E\zeta([0,\tau_1]\times\cdots\times [0,\tau_p])^m =\log  \rho 
 \label{a1.10}
\end{equation} where $\tau_1,\cdots,\tau_p$ are i.i.d. exponential times with
parameter 1 independent of $X$.

In the rest of the paper, we use $\tau_1,\cdots, \tau_p$ to represent
independent exponential times with mean 1, and we use $\Sigma_n$ for the set of all
permutations on $\{1,\cdots, n\}$.  We assume  that
$\{\tau_1,\cdots, \tau_p\}$ and $\{X_1(t),\cdots, X_p(t)\}$ are independent. We begin with a useful  
representation of the $m$'th moment of the random variable
\begin{equation}
\zeta \big([0,\tau_1]\times\cdots\times [0,\tau_p]\big).\label{a2.1}
\end{equation}
 Write $\psi(\lambda)=|\la|^{\bb}$ and $Q(\lambda)=\big[1+\psi(\lambda)\big]^{-1}$.

\begin{lemma}\label{lem-exprep}
\begin{eqnarray} &&
\E\Big[\zeta^{z}([0,\tau_1]\times\cdots\times [0,\tau_p])^m\Big]\label{a2.8}\\ &&  =
\int_{(\R^d)^m}e^{i\sum_{k=1}^m\la_{k}\cdot z}
\bigg[\sum_{\pi\in \Sigma_m}
\prod_{k=1}^m Q\Big(\sum_{j=1}^k\lambda_{\pi(j)}\Big)\bigg]^p
\prod_{k=1}^m\phi_{d-\si}(\lambda_k)\,d\lambda_k\nonumber
\end{eqnarray}
and for any fixed $t_{1},\ldots,t_{p}>0$
\begin{equation}
\E\Big[\zeta ([0,t_1]\times\cdots\times [0,t_p])^n\Big]\leq (t_1\cdots t_p)^ {{\bb p-\si\over\bb p}n}\E\Big[\zeta
\big([0,1]^p\big)^n\Big].\label{a2.jn}
\end{equation}
\end{lemma}

The proof of Lemma \ref{lem-exprep} is given in Section \ref{sec-tidy}.

{\bf  Proof of Theorem \ref{theo-jc} :} Using the multi-parameter version of Kolmogorov's Lemma it suffices to show that we can find  $\de>0$ such that for all $n$ and $M$ we can find a $C<\ff$ such that
\begin{eqnarray}
&&
\E\Big[\Big| \zeta^{z} ([0,t_1]\times\cdots\times [0,t_p])-\zeta^{z'} ([0,t'_1]\times\cdots\times [0,t'_p])\Big |^n\Big]\label{jc.2}\\
&&   \leq C|(z,t_{1},\ldots, t_{p})-(z',t'_{1},\ldots, t'_{p})|^{\de n}\nonumber
\end{eqnarray}
uniformly in $(z,t_{1},\ldots, t_{p}),(z',t'_{1},\ldots, t'_{p})\in R^{d}\times [0,M]^{p}$. To this end it suffices to show separately that 
\begin{eqnarray}
&&
\E\Big[\Big| \zeta^{z} ([0,t_1]\times\cdots\times [0,t_p])-\zeta^{z'} ([0,t_1]\times\cdots\times [0,t_p])\Big |^n\Big]\label{jc.3}\\
&&   \leq C|(z,t_{1},\ldots, t_{p})-(z',t'_{1},\ldots, t'_{p})|^{\de n}\nonumber
\end{eqnarray}
uniformly in $z, z'\in R^{d}, (t_{1},\ldots, t_{p})\in  [0,M]^{p}$ and 
\begin{eqnarray}
&&
\E\Big[\Big| \zeta^{z} ([0,t_1]\times\cdots\times [0,t_p])-\zeta^{z} ([0,t'_1]\times\cdots\times [0,t'_p])\Big |^n\Big]\label{jc.4}\\
&&   \leq C|(z,t_{1},\ldots, t_{p})-(z',t'_{1},\ldots, t'_{p})|^{\de n}\nonumber
\end{eqnarray}
uniformly in $z\in R^{d}, (t_{1},\ldots, t_{p}), (t'_{1},\ldots, t'_{p})\in  [0,M]^{p}$.   

For (\ref{jc.3}) we note first that by the Mean Value Theorem, for any $u,v\geq 0$ we have 
$|u^{-\si}-v^{-\si}|\leq \si |u-v|\max (u^{-\si-1},v^{-\si-1})$. Applying this to $u=|x-z|,\,v=|x-z'|$ we obtain
\begin{equation}
\Big||x-z|^{-\si}- |x-z'|^{-\si}\Big|\leq C|z-z'|
\( |x-z|^{-\si-1}+ |x-z'|^{-\si-1}\).\label{jc.5}
\end{equation}
Interpolating this with the obvious bound
\begin{equation}
\Big||x-z|^{-\si}- |x-z'|^{-\si}\Big|\leq 
\( |x-z|^{-\si}+ |x-z'|^{-\si}\)\label{jc.6}
\end{equation}
we see that for any $0\leq \de\leq 1$
\begin{equation}
\Big||x-z|^{-\si}- |x-z'|^{-\si}\Big|\leq C|z-z'|^{\de}
\( |x-z|^{-\si-\de}+ |x-z'|^{-\si-\de}\).\label{jc.7}
\end{equation}
Then writing
\begin{equation}
\zeta^{z}_{\si}([0,t_1]\times\cdots\times [0,t_p])=\int_0^{t_1}\cdots\int_0^{t_p}
\vert X_1(s_1)+\cdots+ X_p(s_p)-z\vert^{-\sigma}ds_1\cdots ds_p\label{a1.1jc}
\end{equation}
 and   setting $\si'=\si+\de$ for $\de>0$ sufficiently small so that $\si'$ satisfies (\ref{a1.2})
 we see that
\begin{eqnarray} 
&&\quad
\E\Big[\Big| \zeta_{\si}^{z} ([0,t_1]\times\cdots\times [0,t_p])-\zeta_{\si}^{z'} ([0,t_1]\times\cdots\times [0,t_p])\Big |^n\Big]\label{jc.8}\\
&&   \leq C^{n}|z-z'|^{\de n}\sup_{z}
\E\Big[ \zeta_{\si'}^{z} ([0,t_1]\times\cdots\times [0,t_p])^n\Big]\nonumber\\
&&   \leq C^{n}e^{pM}|z-z'|^{\de n}\sup_{z}
\E\Big[ \zeta_{\si'}^{z} ([0,\tau_1]\times\cdots\times [0,\tau_p])^n\Big]\nonumber\\
&&   \leq C^{n}e^{pM}|z-z'|^{\de n}
\int_{(\R^d)^m}
\bigg[\sum_{\pi\in \Sigma_m}
\prod_{k=1}^m Q\Big(\sum_{j=1}^k\lambda_{\pi(j)}\Big)\bigg]^p
\prod_{k=1}^m\phi_{d-\si'}(\lambda_k)\,d\lambda_k\nonumber
\end{eqnarray}
where the last step used (\ref{a2.8}).
By Jensen's inequality,
\begin{eqnarray} &&
\int_{(\R^d)^n}
\Big[\sum_{\sigma\in\Sigma_n}\prod_{k=1}^nQ\Big(\sum_{j=1}^k
\lambda_{\sigma(j)}\Big)\Big]^p\prod_{i=1}^n\phi_{ d-\si'}(\lambda_{ i})d\lambda_i\label{a2.10}\\ && 
\le (n!)^{p-1}
\sum_{\sigma\in\Sigma_n}\int_{(\R^d)^n}
\prod_{k=1}^nQ^p\Big(\sum_{j=1}^k\lambda_{\sigma(j)}\Big)\phi_{
d-\si'}(\lambda_{ i})d\lambda_i\nonumber\\ &&  =(n!)^p\int_{(\R^d)^n}
\prod_{k=1}^nQ^p(\lambda_k)\phi_{
d-\si'}(\lambda_{i}-\lambda_{i-1})d\lambda_i\nonumber\\ &&
\le (n!)^p\bigg(\int_{\R^d}\phi_{
d-\si'}(\lambda)Q^p(\lambda)d\lambda\bigg)^n\nonumber
\end{eqnarray} where the second step follows from variable substitution and we
used the fact  that
\begin{equation}  
\sup_{\lambda' }\int_{\R^d}\phi_{
d-\si'}(\lambda'-\lambda)Q^p(\lambda)d\lambda=\int_{\R^d}\phi_{
d-\si'}(\lambda)Q^p(\lambda)d\lambda.\label{a2.10a}
\end{equation}
This comes from the fact that the convolution of two positive spherically symmetric and monotone decreasing functions has its maximum at the origin. 
(It suffices to prove this for simple functions, and then for  indicator functions of balls centered at the origin in which case it is obvious.)
Finally, (\ref{a2.10a}) is bounded if $\si'$ satisfies (\ref{a1.2}).  This completes the proof of (\ref{jc.3}).

For (\ref{jc.4}) we note first that it suffices to prove a similar bound in which we vary only one of the $t_{j}$. For definiteness we vary $t_{1}$. By H\"{o}lder's inequality, for any positive
function $f$ and any conjugate $r,r'$
\begin{equation}
\int_{A}f(s_{1},\ldots, s_{p})\,ds_{1}\ldots \,ds_{p}
\leq |A|^{1/r'}\(\int_{A}f^{r}(s_{1},\ldots, s_{p})\,ds_{1}\ldots \,ds_{p}\)^{1/r}\label{jc.10}
\end{equation}
where $|A|$ denotes the Lebesgue measure of $A\subseteq R^{p}$. Hence with $t_1>t'_1$
\begin{eqnarray}
&&\quad|\zeta^{z} ([0,t_1]\times\cdots\times [0,t_p])-\zeta^{z} ([0,t'_1]\times [0,t_2]\times \cdots\times [0,t_p])|
\label{jc.11}\\
&& = \zeta^{z} ([t'_1,t_1]\times [0,t_2]\times\cdots\times [0,t_p]) \nonumber\\
&&\leq M |t_1-t'_1  |^{1/r'} 
\(\int_0^{t_1}\cdots\int_0^{t_p}
\vert X_1(s_1)+\cdots+ X_p(s_p)\vert^{-r\sigma}ds_1\cdots ds_p\)^{1/r}. \nonumber
\end{eqnarray}
Choose a rational $r>1$ so that $r\sigma$ satisfies (\ref{a1.2}). Then we can find arbitrarily large $n$ so that $n/r$ is an integer. For such $n$ we can obtain (\ref{jc.4}) as above, and this is enough for  Kolmogorov's Lemma. (In fact, using H\"{o}lder's inequality we can then obtain 
(\ref{jc.4}) for all $n$.)
\qed

We state (\ref{a1.10}) as a  theorem. The  proof is given in Sections \ref{sec-lb} - \ref{sec-Mlim}.

\bt\label{theo-rho} Under (\ref{a1.2}),
\begin{equation}
\lim_{n\to\infty}{1\over n}\log {1\over (n!)^p}
\E\Big[\zeta ([0,\tau_1]\times\cdots\times [0,\tau_p])^n\Big] =\log \rho \label{a2.9}
\end{equation} where $\rho >0$ is given in (\ref{a1.7}).
\et

\bigskip

The hard part of Theorem \ref{theo-rho} is the upper bound. However, it is easy to obtain a rough upper bound using (\ref{a2.10}). Since we will need this in the proof of Theorem \ref{theo-rho} we state this rough upper bound as a Lemma.

\begin{lemma}\label{lem-rub}
\begin{eqnarray} &&\quad
\lim_{n\to\infty}{1\over n}\log {1\over (n!)^p}
\int_{(\R^d)^n}
\Big[\sum_{\sigma\in\Sigma_n}\prod_{k=1}^nQ\Big(\sum_{j=1}^k
\lambda_{\sigma(j)}\Big)\Big]^p\prod_{i=1}^n\phi_{ d-\si}(\lambda_{ i})d\lambda_i \label{a2.9s}\\
&&\hspace{2.5 in}\le\log\bigg(\int_{\R^d}\phi_{
d-\si}(\lambda)Q^p(\lambda)d\lambda\bigg).\nn
\end{eqnarray} 
\end{lemma}

Unfortunately, by examing the argument in (\ref{a2.10})-(\ref{a2.10a}), it is not
hard to see that we do not obtain  the correct constant.
\bigskip

We now show that Theorem \ref{theo-1} follows from Theorem \ref{theo-rho}.

\noindent {\bf Proof of Theorem \ref{theo-1}.} 
Using (\ref{a2.jn})
\begin{eqnarray} &&
\E\Big[\zeta  ([0,\tau_1]\times\cdots\times [0,\tau_p])^n\Big]\label{a2.13}\\ && 
=\int_0^\infty\!\!\cdots\int_0^\infty e^{-(t_1+\cdots +t_p)}
\E\Big[\zeta  ([0,t_1]\times\cdots\times [0,t_p])^n\Big]dt_1\cdots
dt_p\nonumber\\ && 
\le\E\Big[\zeta ([0,1]^p)^n\Big]
\int_0^\infty\!\!\cdots\int_0^\infty  (t_1\cdots t_p)^{{\bb p-\si\over \bb p}n}
e^{-(t_1+\cdots +t_p)}dt_1\cdots dt_p\nonumber\\ &&  =\E\Big[\zeta
([0,1]^p)^n\Big]\Big[\Gamma
\Big({\bb p-\si\over \bb p}n+1\Big)\Big]^p.\nonumber
\end{eqnarray} By Theorem \ref{theo-rho} and Stirling's formula,
\begin{equation}
\liminf_{n\to\infty}{1\over n}\log { 1\over (n!)^{\si/\bb}}
\E\Big[\zeta  ([0,1]^p)^n\Big]\ge
\log \Big({\bb p\over \bb p -\si}\Big)^{\bb p -\si\over\bb} +\log \rho.\label{a2.14}
\end{equation}
\medskip

On the other hand, notice that $\bar{\tau }\equiv\min\{\tau_1,\cdots,\tau_p\}$
has an exponential distribution with the parameter
$p$. Hence,
\begin{eqnarray} &&\qquad
\E\Big[\zeta\big([0,\tau_1]\times\cdots\times  [0,\tau_p]\big)
\Big]^n\ge \E\Big[\zeta\big([0,\bar{\tau} ]^p\big)^n\Big] =\E\bar{\tau}^{{\bb
p-\si\over\bb}n}\E\Big[\zeta\big([0,1]^p\big)^n\Big]\label{a2.15}\\ &&
=p^{-{\bb p-\si\over\bb}n-1}\Gamma\Big(1+{\bb p-\si\over\bb}n\Big)
\E\Big[\zeta\big([0,1]^p\big)^n\Big] \nonumber
\end{eqnarray} where the second step follows from (\ref{a1.3}). By Stirling's formula
we have
\begin{equation}
\limsup_{n\to\infty}{1\over n}\log { 1\over (n!)^{\si/\bb}}
\E\Big[\zeta ([0,1]^p)^n\Big]\le
\log \Big({\bb p\over \bb p -\si}\Big)^{\bb p -\si\over\bb} +\log \rho.\label{a2.15a}
\end{equation}

Combining (\ref{a2.14}) and (\ref{a2.15a}) gives
\begin{equation}
\lim_{n\to\infty}{1\over n}\log (n!)^{-\si/\bb}
\E\Big[\zeta ([0,1]^p)^n\Big]=
\log \Big({\bb p\over \bb p -\si}\Big)^{\bb p -\si\over\bb} +\log \rho.\label{a2.16}
\end{equation}
 Finally, Theorem \ref{theo-1} follows from Lemma
\ref{lem-1.4}.\qed

\section{ Lower bound for Theorem \ref{theo-1}}\label{sec-lb}

  In this section we prove
\begin{equation}
\liminf_{n\to\infty}{1\over n}\log {1\over (n!)^p}
\E\Big[\zeta ([0,\tau_1]\times\cdots\times [0,\tau_p])^n\Big]
\ge\log \rho.\label{a2.17}
\end{equation}

Our starting point is (\ref{a2.8}). Let $q>1$ be the conjugate of $p$
defined by $p^{-1}+q^{-1}=1$ and let $f$ be a  symmetric, continuous, and strictly positive function on $\R^d$ with $\vert\vert f\vert\vert_{ q,\phi_{ d-\si}}=1$, where
\begin{equation}
\vert\vert f\vert\vert_{ q,\phi_{ d-\si}}=\(\int |f(\lambda)|^{ q}\phi_{
d-\si}(\lambda)d\lambda\)^{ 1/q}.\label{a2.17a}
\end{equation} 
We have
\begin{eqnarray} &&
\bigg(\int_{(\R^d)^n}
\Big[\sum_{\sigma\in\Sigma_n}\prod_{k=1}^nQ\Big(\sum_{j=1}^k
\lambda_{\sigma(j)}\Big)\Big]^p\prod_{i=1}^n\phi_{ d-\si}(\lambda_{ i})d\lambda_i
\bigg)^{1/p}\label{a2.18}\\ && 
\ge\int_{(\R^d)^n}
\sum_{\sigma\in\Sigma_n}\prod_{k=1}^nQ\Big(\sum_{j=1}^k
\lambda_{\sigma(j)}\Big)\prod_{i=1}^nf(\lambda_{ i})\phi_{ d-\si}(\lambda_{
i})d\lambda_i\nonumber\\ && 
=n!\int_{(\R^d)^n}
\prod_{k=1}^nQ\Big(\sum_{j=1}^k
\lambda_j\Big)\prod_{i=1}^nf(\lambda_{ i})\phi_{ d-\si}(\lambda_{
i})d\lambda_i\nonumber\\ &&  =n!\int_{(\R^d)^n}
\prod_{k=1}^nf(\lambda_k-\lambda_{k-1})\phi_{d-\si}(\lambda_k-\lambda_{k-1})
Q(\lambda_k)d\lambda_1\cdots
d\lambda_n\nonumber
\end{eqnarray} 
where we follow the convention that $\lambda_0=0$.
\medskip 

Define the linear operator $T$ on ${\cal L}^2(\R^d)$ as
\begin{equation} Tg(\lambda)=\sqrt{Q(\lambda)}\int_{\R^d}f(\gamma
-\lambda)\phi_{d-\si}(\gamma-\lambda)\sqrt{Q(\gamma)}
g(\gamma)d\gamma\hskip.2in g\in{\cal L}^2(\R^d).\label{a2.19}
\end{equation} To show that $T$ is well defined and continuous on ${\cal
L}^2(\R^d)$, we need only to prove that there is a constant $C>0$ such that
\begin{equation}
\langle h, Tg\rangle\le C\vert\vert g\vert\vert_2\vert\vert h\vert\vert_2
\hskip.2in g,h\in{\cal L}^2(\R^d).\label{a2.20}
\end{equation}
But 
\begin{eqnarray} &&
\langle h, Tg\rangle =\int\!\!\int_{\R^d\times\R^d}f(\gamma
-\lambda)\phi_{d-\si}(\gamma-\lambda)
\sqrt{Q(\lambda)}h(\lambda)
\sqrt{Q(\gamma)}g(\gamma)d\lambda d\gamma\label{a2.21}\\ && 
=\int_{\R^d}f(\gamma)\phi_{d-\si}(\gamma)d\gamma\int_{\R^d}\sqrt{Q(\lambda)}h(\lambda)
\sqrt{Q(\lambda + \gamma)}g(\lambda +\gamma)d\lambda\nonumber\\ && 
\le\Bigg\{\int_{\R^d}\phi_{d-\si}(\gamma)\bigg[\int_{\R^d}\sqrt{Q(\lambda)}h(\lambda)
\sqrt{Q(\lambda + \gamma)}g(\lambda 
+\gamma)d\lambda\bigg]^pd\gamma\Bigg\}^{1/p}.\nonumber
\end{eqnarray} Hence by (\ref{sob.1}), 
$\langle h, Tg\rangle\le \vert\vert Q
\vert\vert_{ pd/\si}\vert\vert g\vert\vert_2\vert\vert h\vert\vert_2$.
\medskip

In addition, one can see that $\langle h, Tg\rangle =\langle g, Th\rangle$ for any
$g,h\in{\cal L}^2(\R^d)$. This means that $T$ is self adjoint. We now let $g$ be a 
bounded and locally supported function on $\R^d$ with $\vert\vert
g\vert\vert_2=1$. Then there is $\delta>0$ such that $f, \phi_{d-\si}, Q\ge\delta$ on the support
of $g$. In addition, notice that $Q\le 1$. Thus,
\begin{eqnarray} &&\hspace{.5in}
\int_{(\R^d)^n}
\prod_{k=1}^nf(\lambda_k-\lambda_{k-1})
\phi_{d-\si}(\lambda_k-\lambda_{k-1})Q(\lambda_k)d\lambda_1\cdots d\lambda_n\label{a2.22}\\ && 
\ge\delta^{4}\vert\vert g\vert\vert_\infty^{-2}\int_{(\R^d)^n}
g(\lambda_1)\nn\\
&&
\Big(\prod_{k=2}^n\sqrt{Q(\lambda_{k-1})}
f(\lambda_k-\lambda_{k-1})\phi_{d-\si}(\lambda_k-\lambda_{k-1})\sqrt{Q(\lambda_k)}\Big)g(\lambda_n)d\lambda_1\cdots d\lambda_n\nonumber\\
&&  =\delta^{4}\vert\vert g\vert\vert_\infty^{-2}\langle g,
T^{n-1}g\rangle.\nonumber
\end{eqnarray}

Consider the spectral representation of the self-adjoint operator $T$:
\begin{equation}
\langle g, Tg\rangle =\int_{-\infty}^\infty \theta\mu_g(d\theta)
\end{equation} where $\mu_g(d\theta)$ is a probability measure on $\R$. 
Therefore
\begin{equation}
\langle g, T^{n-1}g\rangle =\int_{-\infty}^\infty \theta^{n-1}\mu_g(d\theta)
\ge \bigg(\int_{-\infty}^\infty \theta\mu_g(d\theta)\bigg)^{n-1} =\langle g,
Tg\rangle^{n-1}\label{a2.23}
\end{equation} where the second step follows from Jensen's inequality.
\medskip

Hence,
\begin{eqnarray} &&
\liminf_{n\to\infty}{1\over n}\log {1\over n!}
\bigg(\int_{(\R^d)^n}\Big[\sum_{\sigma\in\Sigma_n}\prod_{k=1}^nQ\Big(\sum_{j=1}^k
\lambda_{\sigma(j)}\Big)\Big]^p\prod_{i=1}^n\phi_{ d-\si}(\lambda_{ i})d\lambda_i
\bigg)^{1/p}\label{a2.24}\\ && 
\ge \log \langle g, Tg\rangle\nn\\ && 
=\log\int\!\!\int_{\R^d\times\R^d}f(\gamma -\lambda)\phi_{ d-\si}(\gamma
-\lambda)\sqrt{Q(\lambda)}
\sqrt{Q(\gamma)}g(\lambda)g(\gamma)d\lambda d\gamma\nonumber\\ && 
=\log\int_{\R^d}f(\lambda)\phi_{ d-\si}(\lambda )
\bigg[\int_{\R^d}\sqrt{Q(\lambda +\gamma)}\sqrt{Q(\gamma)} g(\lambda
+\gamma)g(\gamma)d\gamma\bigg]d\lambda.\nonumber
\end{eqnarray} 
Notice that the set of all bounded, locally supported $g$ is dense in 
${\cal L}^2(\R^d)$. Taking the supremum over $g$ on the right hand sides gives
\begin{eqnarray} &&\qquad
\liminf_{n\to\infty}{1\over n}\log {1\over n!}
\bigg(\int_{(R^d)^n}
\Big[\sum_{\sigma\in\Sigma_n}\prod_{k=1}^nQ\Big(\sum_{j=1}^k
\lambda_{\sigma(j)}\Big)\Big]^p\prod_{i=1}^n\phi_{ d-\si}(\lambda_{ i})d\lambda_i
\bigg)^{1/p}\label{a2.25}\\ &&
 \ge \log\sup_{\| g\|_2=1}
\int_{\R^d}f(\lambda)\phi_{ d-\si}(\lambda )
\bigg[\int_{\R^d}\sqrt{Q(\lambda +\gamma)}\sqrt{Q(\gamma)} g(\lambda
+\gamma)g(\gamma)d\gamma\bigg]d\lambda.\nonumber
\end{eqnarray}
Notice that for any $g$, the function 
\begin{equation} H(\lambda)=\int_{\R^d}\sqrt{Q(\lambda
+\gamma)}\sqrt{Q(\gamma)} g(\lambda
+\gamma)g(\gamma)d\gamma\label{a2.26}
\end{equation} is symmetric: $H(-\lambda)=H(\lambda)$. Hence, taking the supremum over 
all symmetric, continuous, and strictly positive functions $f$ with $\vert\vert f\vert\vert_{ q,\phi_{
d-\si}}=1$ on the right gives
\begin{eqnarray} &&\qquad
\liminf_{n\to\infty}{1\over n}\log {1\over n!}
\bigg(\int_{(\R^d)^n}
\Big[\sum_{\sigma\in\Sigma_n}\prod_{k=1}^nQ\Big(\sum_{j=1}^k
\lambda_{\sigma(j)}\Big)\Big]^p\prod_{i=1}^n\phi_{ d-\si}(\lambda_{ i})d\lambda_i
\bigg)^{1/p}\label{a2.27}\\ && 
\ge {1\over p}\log\sup_{\| g\|_2=1}\int_{\R^d}\phi_{ d-\si}(\lambda )
\bigg[\int_{\R^d}\sqrt{Q(\lambda +\gamma)}\sqrt{Q(\gamma)} g(\lambda
+\gamma)g(\gamma)d\gamma\bigg]^pd\lambda\nonumber\\ &&  ={1\over
p}\log\rho.\nonumber
\end{eqnarray} From the relation (\ref{a2.8}), we have proved (\ref{a2.17}).\qed

\section{Proof of Lemma \ref{lem-exprep}}\label{sec-tidy}

Before proving the upper bound for Theorem 2.1 we provide the proof of Lemma \ref{lem-exprep},
since we will need several easy generalizations of this proof. In the course of our proof we will use certain Sobolev-type inequalities which are proven in the Appendix.

We first look at 
\begin{eqnarray}
&&
E\(\prod_{j=1}^{n}\int_0^{\tau_1}
{\vert X_1(s_j)+a_{j}\vert^{-\sigma}}ds_j\)\label{v.1}\\
&&= \sum_{\pi\in\Si_{n}}  
E\(\int_{0\leq s_{\pi(1)}\leq\cdots\leq s_{\pi(n)}\leq \tau_1} \prod_{j=1}^{n}
{\vert X_1(s_j)+a_{j}\vert^{-\sigma}}ds_1\cdots ds_n\)\nonumber\\
&&= \sum_{\pi\in\Si_{n}}  
E\(\int_{0\leq s_{\pi(1)}\leq\cdots\leq s_{\pi(n)}\leq \tau_1}\right.\nonumber\\
&&\left.\hspace{.2 in}\(\int \prod_{j=1}^{n}
{\vert x_{\pi(j)}+a_{\pi(j)}\vert^{-\sigma}}p_{s_{\pi(j)}-s_{\pi(j-1)}}( x_{\pi(j)}-x_{\pi(j-1)})\,dx_{j}\)ds_1\cdots ds_n\)\nonumber\\
&&= \sum_{\pi\in\Si_{n}} \int \prod_{j=1}^{n}
{\vert x_{\pi(j)}+a_{\pi(j)}\vert^{-\sigma}}\prod_{j=1}^{n}u^{1}( x_{\pi(j)}-x_{\pi(j-1)})\,dx_{j} \nonumber
\end{eqnarray}
Similarly, proceeding inductively we obtain
\begin{eqnarray}
&&
E\(\prod_{j=1}^{n}\int_0^{\tau_1}\cdots\int_0^{\tau_p}
{\vert X_1(s_{1,j})+\cdots+ X_p(s_{p,j})-z\vert^{-\sigma}}\prod_{l=1}^{p}\,ds_{l,j}\)\label{v.1a}\\
&&= \sum_{\pi_{1},\ldots, \pi_{p}\in\Si_{n}} \int \prod_{j=1}^{n}
{\vert x_{1,\pi_{1}(j)}+\cdots+x_{p,\pi_{p}(j)}-z\vert^{-\sigma}}\nn\\
&& \hspace{2 in}   \prod_{l=1}^{p}\, \prod_{j=1}^{n}u^{1}( x_{l,\pi_{l}(j)}-x_{l,\pi_{l}(j-1)})\,dx_{l,j}. \nonumber
\end{eqnarray}


For  $f\in \mathcal{S}(R^{d})$ let us consider
\begin{eqnarray}
&&\sum_{\pi_{1},\ldots, \pi_{p}\in\Si_{n}} \int \prod_{j=1}^{n}
{\vert x_{1,\pi_{1}(j)}+\cdots+x_{p,\pi_{p}(j)}-z\vert^{-\sigma}}\label{v.1b}\\
&& \hspace{2 in}   \prod_{l=1}^{p}\, \prod_{j=1}^{n}f( x_{l,\pi_{l}(j)}-x_{l,\pi_{l}(j-1)})\,dx_{l,j}. \nonumber
\end{eqnarray}
By  (\ref{a1.4})  
\begin{equation}
\int 
{\vert x\vert^{-\sigma}}f( x)\,dx=\int \phi_{d-\si}(\la)\widehat{f}( \la)\,d\la
\label{v.3}
\end{equation}
and hence 
\begin{eqnarray}
&&
\int 
{\vert x+a\vert^{-\sigma}}f( x)\,dx=\int e^{i\la\cdot a}\phi_{d-\si}(\la)\widehat{f}( \la)\,d\la\label{v.4}\\
&&=\int \(\int e^{i\la\cdot (x+a)}f( x)\,dx\)\phi_{d-\si}(\la)\,d\la.\nn
\end{eqnarray}
Therefore
\begin{eqnarray} 
&&\qquad\sum_{\pi_{1},\ldots, \pi_{p}\in\Si_{n}} \int \prod_{j=1}^{n}
{\vert x_{1,\pi_{1}(j)}+\cdots+x_{p,\pi_{p}(j)}-z\vert^{-\sigma}}\label{v.1c}\\
&& \hspace{2 in}   \prod_{l=1}^{p}\, \prod_{j=1}^{n}f( x_{l,\pi_{l}(j)}-x_{l,\pi_{l}(j-1)})\,dx_{l,j}. \nonumber\\
&&=\sum_{\pi_{1},\ldots, \pi_{p}\in\Si_{n}} \int e^{-i\sum_{j=1}^n\la_{j}\cdot z}
 \(\int  e^{i\sum_{j=1}^{n}\la_{j}\cdot (x_{1,\pi_{1}(j)}+\cdots+x_{p,\pi_{p}(j)})} \right.\nn\\
&& \hspace{1.3 in}\left. \prod_{l=1}^{p}\, \prod_{j=1}^{n}f( x_{l,\pi_{l}(j)}-x_{l,\pi_{l}(j-1)})\,dx_{l,j}\)
\prod_{j=1}^{n}
\phi_{d-\si}(\la_{j})\,d\la_{j}. \nonumber\\
&&=\sum_{\pi_{1},\ldots, \pi_{p}\in\Si_{n}} \int  
\prod_{l=1}^{p}\(\int  e^{i\sum_{j=1}^{n}\la_{j}\cdot x_{l,\pi_{l}(j)}}  \, \prod_{j=1}^{n}f( x_{l,\pi_{l}(j)}-x_{l,\pi_{l}(j-1)})\,dx_{l,j}\)
\nn\\
&& \hspace{3 in}
e^{-i\sum_{j=1}^n\la_{j}\cdot z}\prod_{j=1}^{n}\phi_{d-\si}(\la_{j})\,d\la_{j}. \nonumber
\end{eqnarray}
Note that with $\si=\pi^{-1}$
\begin{eqnarray}
&&\int  e^{i\sum_{j=1}^{n}\la_{j}\cdot x_{\pi(j)}}  \, \prod_{j=1}^{n}f( x_{\pi(j)}-x_{\pi(j-1)})\,dx_{j}
\label{v.1d}\\
&& = \int  e^{i\sum_{j=1}^{n}\la_{\si (j)}\cdot x_{j}}  \, \prod_{j=1}^{n}f( x_{j}-x_{j-1})\,dx_{j} \nonumber\\
&& = \int  e^{i\sum_{j=1}^{n}(\sum_{k=j}^{n}\la_{\si (k)})\cdot x_{j}}  \, \prod_{j=1}^{n}f( x_{j})\,dx_{j} \nonumber\\
&& =  \prod_{j=1}^{n}\widehat{f}\(\sum_{k=j}^{n}\la_{\si (k)}\) =  \prod_{j=1}^{n}\widehat{f}\(\sum_{k=1}^{j}\la_{\si' (k)}\)\nn
\end{eqnarray}
with $\si'$ defined so that $\si'(j)=\si(n-j),\,\forall j$. 
Hence we obtain
\begin{eqnarray}
&&\sum_{\pi_{1},\ldots, \pi_{p}\in\Si_{n}} \int \prod_{j=1}^{n}
{\vert x_{1,\pi_{1}(j)}+\cdots+x_{p,\pi_{p}(j)}-z\vert^{-\sigma}}\label{v.1cc}\\
&& \hspace{1.5 in}   \prod_{l=1}^{p}\, \prod_{j=1}^{n}f( x_{l,\pi_{l}(j)}-x_{l,\pi_{l}(j-1)})\,dx_{l,j} \nonumber\\
&&=\sum_{\pi_{1},\ldots, \pi_{p}\in\Si_{n}} \int e^{-i\sum_{j=1}^n\la_{j}\cdot z}
\prod_{l=1}^{p}\( \prod_{j=1}^{n}\widehat{f}\(\sum_{k=1}^{j}\la_{\pi_{l} (k)}\)\)
 \prod_{j=1}^{n}
\phi_{d-\si}(\la_{j})\,d\la_{j} \nonumber\\
&&= \int e^{-i\sum_{j=1}^n\la_{j}\cdot z}
\Bigg[\sum_{\pi\in\Si_{n}} \prod_{j=1}^{n}\widehat{f}\(\sum_{k=1}^{j}\la_{\pi (k)}\)
\Bigg]^{p}
\prod_{j=1}^{n}\phi_{d-\si}(\la_{j})\,d\la_{j}. \nonumber
\end{eqnarray}
Assuming that $f,\,\widehat{f}\geq 0$ we see as in (\ref{a2.10}) that
\begin{eqnarray}
&& \int 
\Bigg[\sum_{\pi\in\Si_{n}} \prod_{j=1}^{n}\widehat{f}\(\sum_{k=1}^{j}\la_{\pi (k)}\)
\Bigg]^{p}
\prod_{j=1}^{n}\phi_{d-\si}(\la_{j})\,d\la_{j}
\label{v.6}\\
&&  \le (n!)^p\bigg(\int_{\R^d}\phi_{
d-\si}(\lambda)\(\widehat{f}(\lambda)\)^p\,d\lambda\bigg)^n  \nonumber
\end{eqnarray}
and by (\ref{v.1cc}) with $n=1$
\begin{equation} 
 \int_{\R^d}\phi_{
d-\si}(\lambda)\(\widehat{f}(\lambda)\)^p\,d\lambda = \int_{\R^d} { 1\over |x_{1}+\cdots+x_{p}|^{\si} }\prod_{j=1}^{p}f(x_{j})\,dx_{j}.\label{v.11}
\end{equation} 
By (\ref{60.2}) with $\si$ replaced by $d-\si$
\begin{equation}
\int {1 \over |x_{1}+\cdots+x_{p} |^{\si}}\prod_{j=1}^{p} f( x_{j})\,dx_{j}
\leq C^{ p}\|f\|^{p}_{ pd/( pd-\si)}.\label{t.3}
\end{equation}
Now, $u^{ 1}( x)$ is integrable, monotone decreasing in $|x|$  and  asymptotic 
at $x=0$  to  $u^{ 0}( x)=C|x|^{-\max (0,(d-\bb))}$. Hence 
\begin{equation}
\|u^{1}\|_{ pd/( pd-\si)}<\ff\label{v.102a}
\end{equation}
 if $(d-\bb)pd/( pd-\si)<d$
which follows from (\ref{a1.2}). Choose some non-negative $g\in \mathcal{S}(R^{d})$ with $\widehat{g}\geq 0$ and 
 $\int g(x)\,dx=1$. Set $g_{\ep}(x)=\ep^{-d}g(x/\ep)$.
For any sequence $\ep_{r} 
\rar 0$ let  $f_{r}=g_{\ep_{r}}\ast (u^{1}\widehat{g_{\ep_{r}}})\in \mathcal{S}(R^{d})$.  We see that
\begin{equation}
\lim_{r\rar\ff}\|u^{1}-f_{r}\|_{ pd/( pd-\si)}=0\label{v.102b}
\end{equation} and 
$\widehat{f_{r}}=\widehat{g_{\ep_{r}}} (\widehat{u^{1}}\ast g_{\ep_{r}})$  converges pointwise to $\widehat{u^{1}}$. 
In view of (\ref{v.1a}) and (\ref{v.1cc}), to prove (\ref{a2.8})  it suffices to show that
\begin{eqnarray}
&&\lim_{r\rar\ff }\sum_{\pi_{1},\ldots, \pi_{p}\in\Si_{n}} \int \prod_{j=1}^{n}
{\vert x_{1,\pi_{1}(j)}+\cdots+x_{p,\pi_{p}(j)}-z\vert^{-\sigma}}\label{v.31}\\
&& \hspace{1.5 in}   \prod_{l=1}^{p}\, \prod_{j=1}^{n}f_{r}( x_{l,\pi_{l}(j)}-x_{l,\pi_{l}(j-1)})\,dx_{l,j} \nonumber\\
&&=\sum_{\pi_{1},\ldots, \pi_{p}\in\Si_{n}} \int \prod_{j=1}^{n}
{\vert x_{1,\pi_{1}(j)}+\cdots+x_{p,\pi_{p}(j)}-z\vert^{-\sigma}}\nn\\
&& \hspace{1.5 in}   \prod_{l=1}^{p}\, \prod_{j=1}^{n}u^{1}( x_{l,\pi_{l}(j)}-x_{l,\pi_{l}(j-1)})\,dx_{l,j}. \nonumber
\end{eqnarray}
and 
\begin{eqnarray}
&&\lim_{r\rar\ff }\int e^{-i\sum_{j=1}^n\la_{j}\cdot z}
\Bigg[\sum_{\pi\in\Si_{n}} \prod_{j=1}^{n}\widehat{f_{r}}\(\sum_{k=1}^{j}\la_{\pi (k)}\)
\Bigg]^{p}
\prod_{j=1}^{n}\phi_{d-\si}(\la_{j})\,d\la_{j}\label{v.32}\\
&&= \int e^{-i\sum_{j=1}^n\la_{j}\cdot z}
\Bigg[\sum_{\pi\in\Si_{n}} \prod_{j=1}^{n}\widehat{u^{1}}\(\sum_{k=1}^{j}\la_{\pi (k)}\)
\Bigg]^{p}
\prod_{j=1}^{n}\phi_{d-\si}(\la_{j})\,d\la_{j}. \nonumber
\end{eqnarray}

For fixed $\pi_{1},\ldots, \pi_{p}\in\Si_{n}$, the  difference between integral on the the 
right hand side of (\ref{v.31}) and the 
left hand side of (\ref{v.31}) for fixed $r$ 
 is 
 \begin{equation}
\int \prod_{j=1}^{n}
{\vert x_{1,\pi_{1}(j)}+\cdots+x_{p,\pi_{p}(j)}-z\vert^{-\sigma}}F_{r}\prod_{l=1}^{p}\, \prod_{j=1}^{n}\,dx_{l,j}\label{v.31diff}
 \end{equation}
 with 
 \begin{equation}
F_{r}=\prod_{l=1}^{p}\, \prod_{j=1}^{n}u^{1}( x_{l,\pi_{l}(j)}-x_{l,\pi_{l}(j-1)})
-\prod_{l=1}^{p}\, \prod_{j=1}^{n}f_{r}( x_{l,\pi_{l}(j)}-x_{l,\pi_{l}(j-1)}).\label{v.31difg}
 \end{equation}
Writing $A_{(l-1)n+j}=u^{1}( x_{l,\pi_{l}(j)}-x_{l,\pi_{l}(j-1)})$, $B_{(l-1)n+j}=f_{r}( x_{l,\pi_{l}(j)}-x_{l,\pi_{l}(j-1)})$, we can write
 \begin{equation}
 F_{r}=\prod_{s=1}^{np}A_{s}-\prod_{s=1}^{np}B_{r,s}=\sum_{t=1}^{np}
\prod_{s=1}^{t-1}A_{s}\(A_{t}-B_{r,t}\)\prod_{s=t+1}^{np}B_{r,s}.
\label{v.31difh} \end{equation}
It suffices to show that
\begin{eqnarray}
&&
\Bigg|\int \prod_{j=1}^{n}
{\vert x_{1,\pi_{1}(j)}+\cdots+x_{p,\pi_{p}(j)}-z\vert^{-\sigma}}\label{v.31difj}\\
&&\hspace{1 in}
\lc\prod_{s=1}^{t-1}A_{s}\(A_{t}-B_{r,t}\)\prod_{s=t+1}^{np}B_{r,s}\rc\prod_{l=1}^{p}\, \prod_{j=1}^{n}\,dx_{l,j}\Bigg|
\nn
\end{eqnarray}
goes to $0$ as $r\rar\ff$. 
It is easy to see that the product in brackets can be written in the form needed for (\ref{60.2n}). 
More precisely,
\begin{equation}
\prod_{s=1}^{t-1}A_{s}\(A_{t}-B_{r,t}\)\prod_{s=t+1}^{np}B_{r,s}
=\prod_{l=1}^{p}H_{l}\label{v.31difk}
\end{equation}
with 
\begin{equation}
H_{l}=\prod_{j=1}^{n}h_{l,j}( x_{l,\pi_{l}(j)}-x_{l,\pi_{l}(j-1)}) \label{v.31difm}
\end{equation}
where
\[h_{l,j}=\left\{\begin{array}{ll}
u^{1}&\mbox{if $(l-1)n+j<t$}\\
u^{1}-f_{r}&\mbox{if $(l-1)n+j=t$}\\
f_{r}&\mbox{if $(l-1)n+j>t$.}
\end{array}
\right.\]
By (\ref{60.2n}) with $\si$ replaced by $d-\si$, we see that (\ref{v.31difj}) is bounded by
\begin{eqnarray}
&&C\prod_{l=1}^{p}\|\prod_{j=1}^{n}h_{l,j}( x_{l,\pi_{l}(j)}-x_{l,\pi_{l}(j-1)})\|_{ pd/( pd-\si)}
\label{v.31difn}\\
&&=C\prod_{l=1}^{p}\|\prod_{j=1}^{n}h_{l,j}( x_{l,\pi_{l}(j)})\|_{ pd/( pd-\si)}   \nonumber\\
&&=C\prod_{l=1}^{p}\prod_{j=1}^{n}\|h_{l,j}\|_{ pd/( pd-\si)}   \nonumber
\end{eqnarray}
Using (\ref{v.102a}) and (\ref{v.102b}) it is easy to see that this goes to $0$ as $r\rar\ff$, completing the proof of (\ref{v.31}).

Let $\|f\|_{p,\phi_{d-\si}}$ denote the $L^{p}$ norm on 
$R^{dn}$ with respect to the measure $\prod_{j=1}^{n}\phi_{d-\si}(\la_{j})\,d\la_{j}$ so that
\begin{eqnarray}
&&
\int 
\Bigg[\sum_{\pi\in\Si_{n}} \prod_{j=1}^{n}h\(\sum_{k=1}^{j}\la_{\pi (k)}\)
\Bigg]^{p}
\prod_{j=1}^{n}\phi_{d-\si}(\la_{j})\,d\la_{j}\label{v.15h}\\
&&=\|\sum_{\pi\in\Si_{n}} \prod_{j=1}^{n}h\(\sum_{k=1}^{j}\la_{\pi (k)}\)\|^{p}_{p,\phi_{d-\si}}\nn
\end{eqnarray}
Then the absolute value of the  difference between the left hand side of (\ref{v.32}) for fixed $r$ and the 
right hand side of (\ref{v.32}) is bounded by
\[\|\sum_{\pi\in\Si_{n}}\lc  \prod_{j=1}^{n}\widehat{u^{1}}\(\sum_{k=1}^{j}\la_{\pi (k)}\)- \prod_{j=1}^{n}\widehat{f_{r}}\(\sum_{k=1}^{j}\la_{\pi (k)}\)\rc\|_{p,\phi_{d-\si}}^{p}\]
and
\begin{eqnarray}
&&     
\|\sum_{\pi\in\Si_{n}}\lc  \prod_{j=1}^{n}\widehat{u^{1}}\(\sum_{k=1}^{j}\la_{\pi (k)}\)- \prod_{j=1}^{n}\widehat{f_{r}}\(\sum_{k=1}^{j}\la_{\pi (k)}\)\rc\|_{p,\phi_{d-\si}} \nonumber\\
&&  \leq n! 
\|  \prod_{j=1}^{n}\widehat{u^{1}}\(\sum_{k=1}^{j}\la_{k}\)- \prod_{j=1}^{n}\widehat{f_{r}}\(\sum_{k=1}^{j}\la_{k}\)\|_{p,\phi_{d-\si}} \nonumber\\
&&  \leq n! \sum_{m=1}^{n}
\|  \lc\prod_{j=1}^{m-1}\widehat{u^{1}}\(\sum_{k=1}^{j}\la_{k}\)\rc
\Big|(\widehat{u^{1}}-\widehat{f_{r}})\(\sum_{k=1}^{m}\la_{k}\)\Big| 
\lc\prod_{j=m+1}^{n}\widehat{f_{r}}\(\sum_{k=1}^{j}\la_{k}\)\rc\|_{p,\phi_{d-\si}}. \nonumber
\end{eqnarray}
As in (\ref{a2.10})
\begin{eqnarray}
&&
\|  \lc\prod_{j=1}^{m-1}\widehat{u^{1}}\(\sum_{k=1}^{j}\la_{k}\)\rc
\Big|(\widehat{u^{1}}-\widehat{f_{r}})\(\sum_{k=1}^{m}\la_{k}\)\Big| 
\lc\prod_{j=m+1}^{n}\widehat{f_{r}}\(\sum_{k=1}^{j}\la_{k}\)\rc\|^{p}_{p,\phi_{d-\si}}
 \nonumber\\
 &&=\int  \lc\prod_{j=1}^{m-1}\phi_{d-\si}(\la_{j-1}-\la_{j})\(\widehat{u^{1}}(\lambda_{j})\)^p\rc
\phi_{d-\si}(\la_{m-1}-\la_{m})\Big|(\widehat{u^{1}}-\widehat{f_{r}})\(\la_{m}\)\Big|^{p} 
\nonumber\\
&&\hspace{1 in}\lc\prod_{j=m+1}^{n}\phi_{d-\si}(\la_{j-1}-\la_{j})\(\widehat{f_{r}}(\lambda_{j})\)^p\rc
\,d\la_{1}\cdots\,d\la_{n} \nonumber\\
 &&\leq\bigg(\int_{\R^d}\phi_{d-\si}(\lambda)\(\widehat{f_{r}}(\lambda)\)^p\,d\lambda\bigg)^{n-m}\int  \lc\prod_{j=1}^{m-1}\phi_{d-\si}(\la_{j-1}-\la_{j})\(\widehat{u^{1}}(\lambda_{j})\)^p\rc 
\nonumber\\
&&\hspace{1.8 in}\phi_{d-\si}(\la_{m-1}-\la_{m})\Big|(\widehat{u^{1}}-\widehat{f_{r}})\(\la_{m}\)\Big|^{p}
\,d\la_{1}\cdots\,d\la_{m} \label{v.20}
\end{eqnarray}
As in (\ref{v.11})-(\ref{t.3}), $\int_{\R^d}\phi_{d-\si}(\lambda)\(\widehat{f_{r}}(\lambda)\)^p\,d\lambda$ is bounded by $C\|f_{r}\|^{p}_{ pd/( pd-\si)}$ so it remains to show that
\begin{eqnarray}
&&\lim_{r\rar\ff}
\int  \lc\prod_{j=1}^{m-1}\phi_{d-\si}(\la_{j-1}-\la_{j})\(\widehat{u^{1}}(\lambda_{j})\)^p\rc \label{v.21}\\
&&\hspace{1 in}
\phi_{d-\si}(\la_{m-1}-\la_{m})\Big|(\widehat{u^{1}}-\widehat{f_{r}})\(\la_{m}\)\Big|^{p}
\,d\la_{1}\cdots\,d\la_{m}=0.\nn
\end{eqnarray}
We use the uniform integrability of \[G_{r}=:
\prod_{j=1}^{m-1}\(\widehat{u^{1}}(\lambda_{j})\)^p\Big|(\widehat{u^{1}}-\widehat{f_{r}})\(\la_{m}\)\Big|^{p}\] with respect to the measure $d\mu=\prod_{j=1}^{m}\phi_{d-\si}(\la_{j-1}-\la_{j})d\la_{j}$. To see that $G_{r}$ is uniformly integrable it suffices to show that for some $\ep>0$
\begin{eqnarray}
&&\int G_{r}^{1+\ep}\,d\mu=
\int  \lc\prod_{j=1}^{m-1}\phi_{d-\si}(\la_{j-1}-\la_{j})\(\widehat{u^{1}}(\lambda_{j})\)^{p(1+\ep)}\rc \label{v.22}\\
&&\hspace{1 in}
\phi_{d-\si}(\la_{m-1}-\la_{m})\Big|(\widehat{u^{1}}-\widehat{f_{r}})\(\la_{m}\)\Big|^{p(1+\ep)}
\,d\la_{1}\cdots\,d\la_{m}\nn\\&& \leq
\int  \lc\prod_{j=1}^{m-1}\phi_{d-\si}(\la_{j-1}-\la_{j})\(\widehat{u^{1}}(\lambda_{j})\)^{p(1+\ep)}\rc \nn\\
&&\hspace{.1 in}
\phi_{d-\si}(\la_{m-1}-\la_{m})\lc\(\widehat{u^{1}}(\lambda_{j})\)^{p(1+\ep)}+\(\widehat{f_{r}}(\lambda_{j})\)^{p(1+\ep)}\rc 
\,d\la_{1}\cdots\,d\la_{m}\nn
\end{eqnarray}
is bounded uniformly in $r$ and this follows as before. 
Since $\lim_{r\rar\ff}G_{r} =0$ we see that (\ref{v.21}) holds and this establishes (\ref{v.32}).

Let
\begin{equation}
u_{n,t}(y_{1},\ldots,y_{n})=\int_{0\leq s_{1}\leq\cdots\leq s_{n}\leq t}
\prod_{j=1}^{n} p_{s_{j}-s_{j-1}}( y_{j}-y_{j-1})ds_j.\label{v.40}
\end{equation}
To prove (\ref{a2.jn}) we first note as in (\ref{v.1})
\begin{eqnarray}
&&
E\(\prod_{j=1}^{n}\int_0^{t_1}
{\vert X_1(s_j)+a_{j}\vert^{-\sigma}}ds_j\)\label{v.41}\\
&&= \sum_{\pi\in\Si_{n}}  
E\(\int_{0\leq s_{\pi(1)}\leq\cdots\leq s_{\pi(n)}\leq t_1} \prod_{j=1}^{n}
{\vert X_1(s_j)+a_{j}\vert^{-\sigma}}ds_1\cdots ds_n\)\nonumber\\
&&= \sum_{\pi\in\Si_{n}}  
\int_{0\leq s_{\pi(1)}\leq\cdots\leq s_{\pi(n)}\leq t_1}\nonumber\\
&&\hspace{.2 in}\(\int \prod_{j=1}^{n}
{\vert x_{\pi(j)}+a_{j}\vert^{-\sigma}}p_{s_{\pi(j)}-s_{\pi(j-1)}}( x_{\pi(j)}-x_{\pi(j-1)})\,dx_{j}\)ds_1\cdots ds_n\nonumber\\
&&= \sum_{\pi\in\Si_{n}} \int \prod_{j=1}^{n}
{\vert x_{\pi(j)}+a_{j}\vert^{-\sigma}}\,\,\,u_{n,t_{1}}(x_{\pi(1)},\ldots,x_{\pi(n)}) \,dx_{j} \nonumber
\end{eqnarray}
Similarly, proceeding inductively we obtain
\begin{eqnarray}
&&
E\(\prod_{j=1}^{n}\int_0^{t_1}\cdots\int_0^{t_p}
{\vert X_1(s_{1,j})+\cdots+ X_p(s_{p,j})-z\vert^{-\sigma}}\prod_{l=1}^{p}\,ds_{l,j}\)\label{v.41a}\\
&&= \sum_{\pi_{1},\ldots, \pi_{p}\in\Si_{n}} \int \prod_{j=1}^{n}
{\vert x_{1,\pi_{1}(j)}+\cdots+x_{p,\pi_{p}(j)}-z\vert^{-\sigma}}\nn\\
&& \hspace{1 in} \prod_{l=1}^{p} 
u_{n,t_{l}}(x_{\pi_{l}(1)},\ldots,x_{\pi_{l}(n)})  \prod_{l=1}^{p}\, \prod_{j=1}^{n}\,dx_{l,j}. \nonumber
\end{eqnarray}
Then as before we can show that
\begin{eqnarray}
&&\sum_{\pi_{1},\ldots, \pi_{p}\in\Si_{n}} \int \prod_{j=1}^{n}
{\vert x_{1,\pi_{1}(j)}+\cdots+x_{p,\pi_{p}(j)}-z\vert^{-\sigma}}\label{v.41cc}\\
&& \hspace{1 in} \prod_{l=1}^{p} 
u_{n,t_{l}}(x_{\pi_{l}(1)},\ldots,x_{\pi_{l}(n)})  \prod_{l=1}^{p}\, \prod_{j=1}^{n}\,dx_{l,j} \nonumber\\
&&= \int  \prod_{l=1}^{p} e^{-i\sum_{j=1}^n\la_{j}\cdot z}
\Bigg[\sum_{\pi\in\Si_{n}}F_{n,t_{l}}\(\la_{\pi (1)},\ldots,\sum_{k=1}^{n}\la_{\pi (k)}\)
\Bigg]
\prod_{j=1}^{n}\phi_{d-\si}(\la_{j})\,d\la_{j}\nonumber
\end{eqnarray}
where
\begin{eqnarray}
&&F_{n,t}\(\la_{1},\ldots,\la_{n}\)
\label{v.43}\\
&&=  \int e^{i\sum_{j=1}^{n}\la_{j}\cdot y_{j}} \int_{0\leq s_{1}\leq\cdots\leq s_{n}\leq t}
\prod_{j=1}^{n} p_{s_{j}-s_{j-1}}( y_{j})\,ds_j\nonumber\\
&&=   \int_{0\leq s_{1}\leq\cdots\leq s_{n}\leq t}
\prod_{j=1}^{n} e^{-(s_{j}-s_{j-1})\psi(\la_{j})}\,ds_j\nonumber
\end{eqnarray}
which is non-negative. It then follows from the generalized H\"{o}lder's inequality that 
\begin{equation}
\E\Big[\zeta([0,t_1]\times\cdots\times [0,t_p])^m\Big]\leq 
 \prod_{l=1}^{p} \(\E\Big[\zeta([0,t_l]^{p})^m\Big]\)^{1/p}\label{v.44}
\end{equation}
and (\ref{a2.jn}) then follows from the scaling relation (\ref{a1.3}).
\qed 

For future reference we note that (\ref{v.41cc}) and the fact that $F_{n,t}$ is non-negative shows that
\begin{equation}
\sup_{z}\E\Big[\zeta^{z}([0,t_1]\times\cdots\times [0,t_p])^m\Big]\leq \E\Big[\zeta([0,t_1]\times\cdots\times [0,t_p])^m\Big]. \label{domin}
\end{equation}

\section{Upper bound for Theorem \ref{theo-1}}\label{sec-ub}

\medskip In this section we prove
\begin{equation}
\limsup_{n\to\infty}{1\over n}\log {1\over (n!)^p}
\E\Big[\zeta ([0,\tau_1]\times\cdots\times [0,\tau_p])^n\Big]
\le\log \rho.\label{32.1}
\end{equation} 

Define the probability density $h$ on $R^d$ as
\begin{equation}
h(x)=C^{-1}\prod_{j=1}^d\Big({2\sin x_j\over x_j}\Big)^2
\hskip.2in x=(x_1,\cdots, x_d)\in R^d\label{32.4}
\end{equation}
where $C>0$ is the normalizing constant:
$$
C=\int_{\R^d}\prod_{j=1}^d\Big({2\sin x_k\over x_k}\Big)^2dx_1\cdots dx_d.
$$
Clearly, $h$ is symmetric. One can verify that the Fourier transform
$\widehat{h}$ is
$$
\widehat{h}(\lambda)=\int_{\R^d}h(x)e^{i\lambda\cdot x}dx
=C^{-1}(2\pi)^d\Big(1_{[-1,1]^d}\ast 1_{[-1,1]^d}\Big)(\lambda).
$$
In particular, $\widehat{h}$ is non-negative, continuous, with compact
support in the set $[-2,2]^d$,  and
\begin{equation}
\widehat{h}(\lambda)\leq \widehat{h}(0)=1.\label{32.5}
\end{equation}

For each $\epsilon >0$, write
$$
h_\epsilon (x)=\epsilon^{-d}h(\epsilon^{-1}x).\hskip.2in x\in R^d
$$

For some constant $k_{d,\si}$ we have 
\begin{equation}
\int_{R^{d}} {k_{d,\si} \over |s-\la|^{d-\si/2}}\,{k_{d,\si} \over |s|^{d-\si/2}}\,ds={C_{d,\si} \over |\la|^{d-\si}}
=\phi_{d-\si}(\lambda).\label{64.1}
\end{equation}
Let
\begin{equation}
\wp_{\al, \ep} ( \lambda)={k_{d,\si}\widehat{h}(\ep \lambda) \over \al+|\la|^{d-\si/2}}\label{64.2}
\end{equation}
and note that by (\ref{32.5}) and (\ref{64.1})
\begin{equation}
\wp_{\al, \ep}\ast \wp_{\al, \ep}( \lambda)\leq \wp_{\al, 0}\ast \wp_{\al, 0}( \lambda)
\leq \phi_{d-\si}( \lambda).\label{64.2a}
\end{equation}
Let
\begin{equation}
\th_{ \al,\ep}( x)=\int e^{ ix\cdot \la}{k_{d,\si}\widehat{h}(\ep \lambda) \over \al+|\la|^{d-\si/2}}\,d\la
=\int e^{ ix\cdot \la}\wp_{\al, \ep} ( \lambda)\,d\la.\label{64.3}
\end{equation}
Then
\begin{equation}
\th^{2}_{ \al,\ep}( x)=\int e^{ ix\cdot \la}\wp_{\al, \ep}\ast \wp_{\al, \ep} ( \lambda)\,d\la.\label{64.4}
\end{equation}

\medskip

 Define
\begin{equation}
\zeta_{ \al,\ep}\big([0,t_1]\times\cdots\times [0,t_p]\big)
=\int_0^{t_1}\cdots\int_0^{t_p}
\th^{2}_{ \al,\ep}(X_1(s_1)+\cdots+ X_p(s_p))ds_1\cdots ds_p.\label{a4.2}
\end{equation}

(\ref{32.1}) will follow from the next two Lemmas.

\begin{lemma}\label{lem-expapprox-al}
\begin{equation}
\limsup_{\al,\ep\to 0^+}\limsup_{n\to\infty}{1\over n}\log{1\over (n!)^p}
\E\Big[(\zeta  -\zeta_{ \al,\ep})\big([0,\tau_1]\times\cdots\times 
[0,\tau_p]\big)
\Big]^n=-\infty.\label{a4.9}
\end{equation}
\end{lemma}

\begin{lemma}\label{lem-al}
\begin{equation}
\limsup_{n\to\infty}{1\over n}\log {1\over (n!)^p}
\E\Big[\zeta_{\al,\ep} ([0,\tau_1]\times\cdots\times [0,\tau_p])^n\Big]
\le\log  \rho .\label{32.2}
\end{equation}
\end{lemma}

{\bf  Proof of Lemma \ref{lem-expapprox-al}}

By (\ref{64.4}) and (\ref{a4.2})
\begin{eqnarray} &&\hspace{.2in}
\zeta_{ \al,\ep}\big([0,t_1]\times\cdots\times [0,t_p]\big)\label{a4.3}\\ &&  =\int_{R^d}\( 
\int_0^{t_1}\!\!\cdots\!\!\int_0^{t_p}
\exp\Big\{i\lambda\cdot \big(X_1(s_1)+\cdots +X_p(s_p)\big)\Big\} ds_1\cdots
ds_p\)\nn\\
&&\hspace{3.5 in} \wp_{\al, \ep}\ast \wp_{\al, \ep}(\lambda)\,d\lambda\nonumber
\end{eqnarray} 
 Following the same procedure used for (\ref{a2.8}),
\begin{eqnarray} &&\qquad
\E\Big[(\zeta -\zeta_{ \al,\ep})\big([0,\tau_1]\times\cdots\times 
[0,\tau_p]\big)^n
\Big]\label{a4.5}\\ &&  =\int_{(\R^d)^n}
\Big[\sum_{\sigma\in\Sigma_n}\prod_{k=1}^nQ\Big(\sum_{j=1}^k
\lambda_{\sigma(j)}\Big)\Big]^p 
\prod_{k=1}^n\big[\phi_{d-\si}(\lambda_k)-\wp_{\al, \ep}\ast \wp_{\al, \ep}(\lambda_k)\big]
\,d\lambda_k\nonumber
\end{eqnarray} 
where $Q(\lambda)=\big[1+\psi(\lambda)\big]^{-1}$.
\medskip

Note  that
\begin{eqnarray}&&
0\leq \phi_{d-\si}(\lambda)-\wp_{\al, \ep}\ast \wp_{\al, \ep}(\lambda)=
\label{a4.7}\\
&&=\(\phi_{d-\si}(\lambda)-\wp_{\al,0}\ast \wp_{\al, 0}(\lambda)\)
+\(\wp_{\al,0}\ast \wp_{\al, 0}(\lambda)-\wp_{\al, \ep}\ast \wp_{\al, \ep}(\lambda)\)\nonumber
\end{eqnarray}
By (\ref{64.1}) we have
\begin{eqnarray}
&&\hspace{.3 in}0\leq \phi_{d-\si}(\lambda)-\wp_{\al,0}\ast \wp_{\al, 0}(\lambda)
\label{64.6}\\
&& =  k_{d,\si}^{2}\(\int {1 \over |s-\la|^{d-\si/2}}\,{1 \over |s|^{d-\si/2}}\,ds
-\int {1 \over \al+ |s-\la|^{d-\si/2}}\,{1 \over \al+|s|^{d-\si/2}}\,ds\)\nonumber\\
&& \leq   C{\al^{\de} \over |\la|^{d-\si/2+\de}}.\nonumber
\end{eqnarray}
We also have
\begin{eqnarray}
&&0\leq \wp_{\al,0}\ast \wp_{\al, 0}(\lambda)-\wp_{\al, \ep}\ast \wp_{\al, \ep}(\lambda)
\label{64.7}\\
&&  =k_{d,\si}^{2}\(\int {1 \over \al+|s-\la|^{d-\si/2}}\,{1 \over \al+|s|^{d-\si/2}}\,ds\right.\nn\\
&&\left.\hspace{1 in}
-\int {\widehat{h}(\ep (s-\la)) \over \al+ |s-\la|^{d-\si/2}}\,{\widehat{h}(\ep s) \over \al+|s|^{d-\si/2}}\,ds\) \nonumber\\
&& \leq k_{d,\si}^{2}\(\int {1 \over |s-\la|^{d-\si/2}}\,{1-\widehat{h}(\ep s) \over |s|^{d-\si/2}}\,ds\right.\nn\\
&&\left.\hspace{1 in}
+\int {1-\widehat{h}(\ep (s-\la)) \over  |s-\la|^{d-\si/2}}\,{1 \over |s|^{d-\si/2}}\,ds\) \nonumber
\end{eqnarray}
Fix $\ga>0$ and choose  $\tau>0$  so that (see (\ref{32.5}) )
\begin{equation}
0\leq (1-\wh{h}(z))\leq \ga,\hspace{.2 in}|z|\leq \tau.\label{32.9}
\end{equation}
By considering separately the regions $s\leq \tau/\ep $ and $s> \tau/\ep $ we see that
\begin{eqnarray} &&
{1-\widehat{h}(\ep s) \over |s|^{d-\si/2}}\leq \ga {1 \over |s|^{d-\si/2}}+\({\ep \over \tau}\)^{\de}
{1 \over |s|^{d-\si/2-\de}}\label{32.10}\\
&&\leq \ga\({1 \over |s|^{d-\si/2}}+   {1 \over |s|^{d-\si/2-\de}}\)\nn
\end{eqnarray}
for $\ep>0$ sufficiently small.    Here we can take any $\de$ sufficiently small with
$\si+\de<\min (d,p\bb) $.
Our Lemma then  follows from Lemma \ref{lem-rub} by  first taking $\al, \ep\rar 0$ with $\ga>0$ fixed  and then letting $\ga\rar 0$.

Our Lemma then  follows from Lemma \ref{lem-rub} by taking $\de>0$ sufficiently small.
  \qed

{\bf  Proof of Lemma \ref{lem-al} }  Define
\begin{equation}
\zeta_{ \al,\ep',\ep}\big([0,t_1]\times\cdots\times [0,t_p]\big)
=\int_0^{t_1}\cdots\int_0^{t_p}
\th_{ \al,\ep'}^{2}\ast h_{\ep}(X_1(s_1)+\cdots+ X_p(s_p))ds_1\cdots ds_p.\label{32.11m}
\end{equation}
Lemma \ref{lem-al} will follow from the next two Lemmas.

\begin{lemma}\label{lem-expapprox-alm}
\begin{equation}
\limsup_{\ep\to 0^+}\limsup_{n\to\infty}{1\over n}\log{1\over (n!)^p}
\E\Big[(\zeta_{ \al,\ep',\ep}  -\zeta_{ \al,\ep'})\big([0,\tau_1]\times\cdots\times 
[0,\tau_p]\big)
\Big]^n=-\infty.\label{64.9}
\end{equation}
\end{lemma}

\begin{lemma}\label{lem-alm}
\begin{equation}
\limsup_{n\to\infty}{1\over n}\log {1\over (n!)^p}
\E\Big[\zeta_{ \al,\ep',\ep} ([0,\tau_1]\times\cdots\times [0,\tau_p])^n\Big]
\le\log  \rho .\label{64.10}
\end{equation}
\end{lemma}

{\bf  Proof of Lemma \ref{lem-expapprox-alm} }  Following the same procedure used for (\ref{a2.8}),
\begin{eqnarray} &&\qquad
\E\Big[(\zeta_{ \al,\ep'}  -\zeta_{ \al,\ep',\ep})\big([0,\tau_1]\times\cdots\times 
[0,\tau_p]\big)^n
\Big]\label{a4.5}\\ &&  =\int_{(\R^d)^n}
\Big[\sum_{\sigma\in\Sigma_n}\prod_{k=1}^nQ\Big(\sum_{j=1}^k
\lambda_{\sigma(j)}\Big)\Big]^p 
\prod_{k=1}^n\big[(1-\widehat{h}(\ep \la))\wp_{\al, \ep}\ast \wp_{\al, \ep}(\lambda_k)\big]
\,d\lambda_k\nonumber\\ && \leq \int_{(\R^d)^n}
\Big[\sum_{\sigma\in\Sigma_n}\prod_{k=1}^nQ\Big(\sum_{j=1}^k
\lambda_{\sigma(j)}\Big)\Big]^p 
\prod_{k=1}^n\big[(1-\widehat{h}(\ep \la))\phi_{d-\si}(\lambda_k)\big]
\,d\lambda_k\nonumber
\end{eqnarray} 
by (\ref{64.2a}) and the proof follows as in the proof of Lemma \ref{lem-expapprox-al}.\qed

{\bf  Proof of Lemma \ref{lem-alm} }

 Define

\begin{equation}
\zeta_{ \al,\ep}^{z}\big([0,t_1]\times\cdots\times [0,t_p]\big)
=\int_0^{t_1}\cdots\int_0^{t_p}
\th_{ \al,\ep}^{2}(X_1(s_1)+\cdots+ X_p(s_p)-z)ds_1\cdots ds_p.\label{32.11}
\end{equation}
 Let $M>0$ be
fixed but arbitrary. By definition, using the fact that both $h_{\ep}(z)$ and $\zeta_{\al,\ep'}^{z}$ are non-negative functions
\begin{eqnarray} &&
\zeta_{\al,\ep',\ep}\big([0,t_1]\times\cdots\times  [0,t_p]\big)\label{b4.10}\\
&&=\sum_{y\in\Z^d}\int_{[0,M]^d}h_{\ep} (yM+z)
\zeta_{\al,\ep'}^{yM+z}\big([0,t_1]\times\cdots\times [0,t_p]\big)dz \nonumber
\\ &&\le\int_{[0,M]^d}\widetilde{h}_{\ep} (z)
\widetilde{\zeta}_{\al,\ep'}^z\big([0,t_1]\times\cdots\times [0,t_p]\big)dz \nonumber
\end{eqnarray} where
\begin{equation}
\widetilde{h}_{\ep}(x)=\sum_{y\in\Z^d}h_{\ep} (yM+z),
\hskip.2in
\widetilde{\zeta}_{\al,\ep'}^z\big([0,t_1]\times\cdots\times [0,t_p]\big)
=\sum_{y\in\Z^d}\zeta_{\al,\ep'}^{yM+z}\big([0,t_1]
\times\cdots\times [0,t_p]\big)\label{b4.11}
\end{equation} are two periodic functions on $\R^d$ with the period $M>0$.

\medskip

By Parseval's identity
\begin{eqnarray} &&
\int_{[0,M]^d}\widetilde{h}_{\ep} (z)
\widetilde{\zeta}_{\al,\ep'}^z\big([0,t_1]\times\cdots\times [0,t_p]\big)dz\label{b4.12}\\ &&
={1\over M^d}\sum_{y\in\Z^d}\bigg(\int_{[0,M]^d}\widetilde{h}_{\ep} (x)
\exp\Big\{-i{2\pi\over M}(y\cdot x)\Big\}dx\bigg) \nonumber
\\ && \hskip.3in\times\bigg(\int_{[0,M]^d}
\widetilde{\zeta}_{\al,\ep'}^x\big([0,t_1]\times\cdots\times [0,t_p]\big)
\exp\Big\{i{2\pi\over M}(y\cdot x)\Big\}dx\bigg).\nonumber
\end{eqnarray} 
By periodicity
\begin{eqnarray} &&
\int_{[0,M]^d}\widetilde{h}_{\ep} (x)
\exp\Big\{-i{2\pi\over M}(y\cdot x)\Big\}dx\label{b4.13}\\ &&
 =\sum_{z\in\Z^d}\int_{[0,M]^d}h_{\ep} (zM+x)
\exp\Big\{-i{2\pi\over M}(y\cdot x)\Big\}dx\nonumber\\ &&
 =\sum_{z\in\Z^d}\int_{zM+[0,M]^d}h_{\ep} (x)
\exp\Big\{-i{2\pi\over M}\big(y\cdot (x-zM\big)\big)\Big\}dx\nonumber\\ &&
 =\sum_{z\in\Z^d}\int_{zM+[0,M]^d}h_{\ep} (x)
\exp\Big\{-i{2\pi\over M}(y\cdot x)\Big\}dx\nonumber\\ &&
 =\int_{\R^d}h_{\ep} (x)
\exp\Big\{-i{2\pi\over M}(y\cdot x)\Big\}dx
=\widehat{h} \Big(\ep {2\pi\over M}y\Big).\nonumber
\end{eqnarray} 
 Similarly, using (\ref{32.11})
\begin{eqnarray} &&
\int_{[0,M]^d}
\widetilde{\zeta}_{\al,\ep'}^x\big([0,t_1]\times\cdots\times [0,t_p]\big)
\exp\Big\{i{2\pi\over M}(y\cdot x)\Big\}dx\label{b4.14}\\ &&=\int_{\R^d}
\zeta_{\al,\ep'}^x\big([0,t_1]\times\cdots\times [0,t_p]\big)
\exp\Big\{i{2\pi\over M}(y\cdot x)\Big\}dx \nonumber
\\ &&=\int_{[0,t_1]\times\cdots\times[0,t_p]}\wp_{\al, \ep'}\ast \wp_{\al, \ep'}( {2\pi\over M}y)\nonumber
\\ &&\hspace{1 in}
\exp\Big\{i{2\pi\over M}y\cdot \big(X_1(s_1)+\cdots +X_p(s_p)\big)
\Big\}ds_1\cdots ds_p. \nonumber
\end{eqnarray}
 Hence,
\begin{eqnarray} &&
\int_{[0,M]^d}\widetilde{h}_{\ep} (z)
\widetilde{\zeta}_{\al,\ep'}^z\big([0,t_1]\times\cdots\times [0,t_p]\big)dz\label{b4.15}\\ &&
={1\over M^d}\sum_{y\in\Z^d}\widehat{h} \Big(\ep {2\pi\over M}y\Big)\wp_{\al, \ep'}\ast \wp_{\al, \ep'}
\Big({2\pi\over M}y\Big)\nn\\
&&
\int_{[0,t_1]\times\cdots\times[0,t_p]}
\exp\Big\{i{2\pi\over M}y\cdot \big(X_1(s_1)+\cdots +X_p(s_p)\big)
\Big\}ds_1\cdots ds_p. \nonumber
\end{eqnarray}
Using the same procedure as the one used to derive Lemma \ref{lem-exprep}, (in fact, here we can proceed more directly, as in \cite{Cadditive}) we can show that
\begin{eqnarray} &&\E\bigg[\int_{[0,M]^d}\widetilde{h}_{\ep} (z)
\widetilde{\zeta}_{\al,\ep'}^z\big([0,\tau_1]\times\cdots\times [0,\tau_p]
\big)dz\bigg]^n\label{b4.16}\\ &&={1\over
M^{dn}}\sum_{y_1,\cdots,y_n\in\Z^d}\bigg(\prod_{k=1}^n\widehat{h} \Big(\ep {2\pi\over M}y_k\Big)\wp_{\al, \ep'}\ast \wp_{\al, \ep'}
\Big({2\pi\over M}y_k\Big)\bigg)\nn\\
&&\hspace{1 in}\bigg[\sum_{\sigma\in\Sigma_n}
\prod_{k=1}^nQ\Big({2\pi\over M}\sum_{j=1}^ky_{\sigma(j)}\Big)\bigg]^p
\nonumber
\end{eqnarray}

By \cite[Theorem 4.1]{Cadditive}, (\ref{32.5}),  (\ref{64.2a}) and the fact that $\widehat{h}$ is
supported in the set $[-2,2]^d$,
\begin{eqnarray} &&
\lim_{n\to\infty}{1\over n}\log {1\over (n!)^p}
\E\bigg[\int_{[0,M]^d}\widetilde{h}_{\ep} (z)
\widetilde{\zeta}_{\al,\ep'}^z\big([0,\tau_1]\times\cdots\times [0,\tau_p]
\big)dz\bigg]^n\label{b4.17}\\ && =\log\Bigg({1\over M^d}\sup_{\|
f\|_{2,\Z^d}=1}\sum_{x\in\Z^d}
\widehat{h} \Big(\ep {2\pi\over M}x\Big)\wp_{\al, \ep'}\ast \wp_{\al, \ep'}
\Big({2\pi\over M}x\Big)\nn\\
&&\hspace{ 1in}
\bigg[\sum_{y\in\Z^d}\sqrt{Q\Big({2\pi\over M}(x+y)\Big)}
\sqrt{Q\Big({2\pi\over M}y\Big)}f(x+y)f(y)\bigg]^p\Bigg)\nonumber
\\ &&\le \log\Big(M^{-d}\rho_M\Big) \nonumber
\end{eqnarray} 
where, setting $a=2\sqrt{d}/\ep$
\begin{eqnarray} &&
\rho_M =\sup_{\| f\|_{2,\Z^d}=1}\,\,\sum_{\vert x\vert\le (2\pi)^{-1}Ma}\wp_{\al,0}\ast \wp_{\al, 0}\( {2\pi\over M}x\)
\label{b4.18}\\
&&\hspace{1 in}
\bigg[\sum_{y\in\Z^d}\sqrt{Q\Big({2\pi\over M}(x+y)\Big)}
\sqrt{Q\Big({2\pi\over M}y\Big)}f(x+y)f(y)\bigg]^p.\nn
\end{eqnarray}
In view of (\ref{b4.10}),
\begin{equation}
\limsup_{n\to\infty}{1\over n}\log {1\over (n!)^p}
\E\Big[\zeta_{\al,\ep',\ep}\big([0,\tau_1]\times\cdots\times  [0,\tau_p]\big)
\Big]^n\le \log\Big(M^{-d}\rho_M\Big).\label{b4.199}
\end{equation}
 By Theorem \ref{theo-alt} of the next section,  letting $M\to\infty$ on the
right hand side gives
\begin{equation}
\limsup_{n\to\infty}{1\over n}\log {1\over (n!)^p}
\E\Big[\zeta_{\al,\ep',\ep}\big([0,\tau_1]\times\cdots\times  [0,\tau_p]\big)
\Big]^n\le \log \rho .\label{b4.20}
\end{equation} 
\qed

\section{ The limit as $M\to\ff$}\label{sec-Mlim} 
\medskip

\bt\label{theo-alt} Let $\rho$ be defined in (1.4) and
$\rho_M$ be defined in (5.10). We have
\begin{equation}
\limsup_{M\to\infty}M^{-d}\rho_M\le \rho.\label{ap.0}
\end{equation}
\et

\proof 
For any $x=(x_1,\cdots, x_d)\in\R^d$, we write
$[x]=([x_1],\cdots, [x_d])$ for the lattice part of $x$ (We also use the the notation
$[\cdots ]$ for parentheses without causing any confusion). For any $f\in {\cal
L}^2(\Z^d)$ with $\| f\|_2=1$,
\begin{eqnarray} &&
\sum_{\vert x\vert\le (2\pi)^{-1}Ma}\wp_{\al,0}\ast \wp_{\al, 0}\({2\pi\over M}x\)
\bigg[\sum_{y\in\Z^d}\sqrt{Q\Big({2\pi\over M}(x+y)\Big)}
\sqrt{Q\Big({2\pi\over M}y\Big)}f(x+y)f(y)\bigg]^p\nn\\ &&
 =\int_{\{\vert \lambda\vert\le (2\pi)^{-1}Ma\}}\wp_{\al,0}\ast \wp_{\al, 0}\({2\pi\over
M}[\lambda]\)
\nn\\
&&\hspace{ 1in}
\bigg[\int_{\R^d}\sqrt{Q\Big({2\pi\over M}([\lambda]+[\gamma])\Big)}
\sqrt{Q\Big({2\pi\over M}[\gamma]\Big)}f([\lambda]+[\gamma])
f([\gamma])d\gamma\bigg]^pd\lambda\nonumber\\ && =\Big({M\over
2\pi}\Big)^d\int_{\{\vert\lambda\vert\le a\}}
\wp_{\al,0}\ast \wp_{\al, 0}\({2\pi\over M}[{M\over 2\pi}\lambda]\)\nonumber\\ &&
\hspace{1 in}\bigg[\Big({M\over 2\pi}\Big)^d
\int_{\R^d}\sqrt{Q_M
\Big(\gamma +{2\pi\over M}\Big[{M\over 2\pi}\lambda\Big]\Big)}
\sqrt{Q_M(\gamma)}  \nonumber\\ &&\hspace{2 in}\times f\Big(\Big[{M\over 2\pi}\lambda\Big]+
\Big[{M\over 2\pi}\gamma\Big]\Big) f\Big(\Big[{M\over
2\pi}\gamma\Big]\Big)d\gamma\bigg]^pd\lambda\label{ap.14}
\end{eqnarray}
 where
\begin{equation} Q_M(\lambda)=Q\Big({2\pi\over M}\Big[{M\over
2\pi}\lambda\Big]\Big)
\hskip.2in\lambda\in\R^d.\label{64.15}
\end{equation}

Write
\begin{equation} g_0(\lambda) =\Big({M\over 2\pi}\Big)^{d/2}f\Big(\Big[{M\over
2\pi}\lambda\Big]\Big)
\hskip.2in\lambda\in\R^d.\label{ap.15}
\end{equation} We have
\begin{equation}
\int_{\R^d}g_0^2(\lambda)d\lambda=\Big({M\over 2\pi}\Big)^d
\int_{\R^d}f^2\Big(\Big[{M\over 2\pi}\lambda\Big]\Big)d\lambda
=\int_{\R^d}f^2([\lambda])d\lambda=\sum_{x\in\Z^d}f^2(x)=1.\label{}
\end{equation} We can also see that under this correspondence,
\begin{equation}
\Big({M\over 2\pi}\Big)^{d/2}f\Big(\Big[{M\over 2\pi}\lambda\Big]+
\Big[{M\over 2\pi}\gamma\Big]\Big)=g_0\Big(\gamma +{2\pi\over M}
\Big[{M\over
2\pi}\lambda\Big]\Big)\hskip.2in\lambda,\gamma\in\R^d.\label{ap.16}
\end{equation}
 
Therefore, we need only to show that for any fixed $a>0$
\begin{eqnarray} &&\qquad
\limsup_{M\to\infty}\sup_{\vert\vert g\vert\vert_2=1}
\int_{\{\vert \lambda\vert\le a\}}\wp_{\al,0}\ast \wp_{\al, 0}\({2\pi\over M}[{M\over
2\pi}\lambda]\)
\label{ap.17}\\ &&
\bigg[\int_{\R^d}\sqrt{Q_M
\Big(\gamma +{2\pi\over M}\Big[{M\over 2\pi}\lambda\Big]\Big)}
\sqrt{Q_M(\gamma)} g\Big(\gamma +{2\pi\over M}
\Big[{M\over 2\pi}\lambda\Big]\Big) g(\gamma)d\gamma\bigg]^pd\lambda
\nonumber\\ &&
 \le\sup_{\vert\vert g\vert\vert_2=1}\int_{\{\vert \lambda\vert\le a\}}
\phi_{d-\si}( \la) 
\bigg[\int_{\R^d}\sqrt{Q(\lambda +\gamma)}
\sqrt{Q(\gamma)}g(\lambda +\gamma)g(\gamma)d\gamma\bigg]^p\,d\lambda.\nonumber
\end{eqnarray} 
To this end, note that by the inverse Fourier transformation  the function
\begin{equation} U_M(\lambda)=\int_{\R^d}\sqrt{Q_M (\gamma +\lambda)}
\sqrt{Q_M(\gamma)}g(\gamma +\lambda) g(\gamma)d\gamma\label{}
\end{equation} is the Fourier transform of the function
\begin{eqnarray} && V_M(x)={1\over
(2\pi)^d}\int_{\R^d}U_M(\lambda)e^{-i\lambda\cdot x}d\lambda\label{ap.18}\\
&&  ={1\over (2\pi)^d}\int_{\R^d}e^{-i\lambda\cdot
x}d\lambda\int_{\R^d}\sqrt{Q_M (\gamma +\lambda)}
\sqrt{Q_M(\gamma)}g(\gamma +\lambda) g(\gamma)d\gamma\nonumber\\ && 
={1\over (2\pi)^d}\int\!\!\int_{\R^d\times\R^d}e^{-i(\lambda-\gamma)\cdot x}
\sqrt{Q_M(\lambda)}g(\lambda)\sqrt{Q_M(\gamma)}g(\gamma)d\lambda
d\gamma\nonumber\\ &&  ={1\over
(2\pi)^d}\bigg\vert\int_{\R^d}e^{ix\cdot\gamma}\sqrt{Q_M(\gamma)}
g(\gamma)d\gamma\bigg\vert^2.\nonumber
\end{eqnarray} Therefore
\begin{eqnarray} &&
\int_{\R^d}\sqrt{Q_M
\Big(\gamma +{2\pi\over M}\Big[{M\over 2\pi}\lambda\Big]\Big)}
\sqrt{Q_M(\gamma)}g\Big(\gamma +{2\pi\over M}
\Big[{M\over 2\pi}\lambda\Big]\Big)g(\gamma)d\gamma \nn\\
&& =U_M\Big({2\pi\over M}\Big[{M\over 2\pi}\lambda\Big]\Big)\nn\\ &&
 ={1\over (2\pi)^d}\int_{\R^d}
\exp\Big\{ix\cdot {2\pi\over M}\Big[{M\over 2\pi}\lambda\Big]\Big\}
\bigg\vert\int_{\R^d}e^{ix\cdot\gamma}\sqrt{Q_M(\gamma)}
g(\gamma)d\gamma\bigg\vert^2dx\nonumber\\ &&
 \le {1\over (2\pi)^d}\int_{\R^d}\Big\vert 1-\exp\Big\{ix\cdot \Big(\lambda -
{2\pi\over M}\Big[{M\over 2\pi}\lambda\Big]\Big)\Big\}\Big\vert\cdot
\bigg\vert\int_{\R^d}e^{ix\cdot\gamma}\sqrt{Q_M(\gamma)}
g(\gamma)d\gamma\bigg\vert^2dx\nonumber\\ &&\hspace{ 1in}
 +{1\over (2\pi)^d}\int_{\R^d} e^{ix\cdot \lambda}
\bigg\vert\int_{\R^d}e^{ix\cdot\gamma}\sqrt{Q_M(\gamma)}
g(\gamma)d\gamma\bigg\vert^2dx.\label{ap.19}
\end{eqnarray} 
By Parseval's identity and by the fact $Q_M\le 1$,
\begin{eqnarray} && {1\over (2\pi)^d}\int_{\R^d}\bigg\vert\int_{\R^d}
e^{ix\cdot\gamma}\sqrt{Q_M(\gamma)}
g(\gamma)d\lambda\bigg\vert^2dx\label{}\\ &&
=\int_{\R^d}Q_M(\gamma)g^2(\gamma)d\gamma
\le\int_{\R^d}g^2(\gamma)d\gamma =1.\nonumber
\end{eqnarray} Hence, the first term on the right hand side of (\ref{ap.19}) tends to
0 uniformly over
$\lambda\in\R^d$ and over all $g\in{\cal L}^2(\R^d)$  with $\vert\vert
g\vert\vert_2=1$ as $M\to\infty$.  The second term on the right hand side of
(\ref{ap.19}) is equal to
\begin{equation}
\int_{\R^d}e^{ix\cdot\lambda}V_M(x)dx=U_M(\lambda)=
\int_{\R^d}\sqrt{Q_M(\lambda +\gamma)}
\sqrt{Q_M(\gamma)}g(\lambda +\gamma)g(\gamma)d\gamma.\label{ap.21}
\end{equation}

Consequently, we will have (\ref{ap.17}) if we can prove
\begin{eqnarray} &&\quad
\limsup_{M\to\infty}\sup_{\vert\vert g\vert\vert_2=1}
\int_{\{\vert \lambda\vert\le a\}}\wp_{\al,0}\ast \wp_{\al, 0}({2\pi\over M}[{M\over
2\pi}\lambda]) \label{ap.22}\\
&&\hspace{ 1in}\times
\bigg[\int_{\R^d}\sqrt{Q_M (\lambda +\gamma)}
\sqrt{Q_M(\gamma)}g(\lambda +\gamma)
g(\gamma)d\gamma\bigg]^pd\lambda\nn\\ &&\le\sup_{\vert\vert
g\vert\vert_2=1}\int_{\{\vert \lambda\vert\le a\}}\phi_{d-\si}(\lambda ) 
\bigg[\int_{\R^d}\sqrt{Q(\lambda +\gamma)}
\sqrt{Q(\gamma)}g(\lambda +\gamma)g(\gamma)d\gamma\bigg]^pd\lambda
\nonumber
\end{eqnarray}

By uniform continuity of the function $Q$ we have that
$Q_M(\cdot)\to Q(\cdot)$ uniformly on $\R^d$. Thus, given $\epsilon >0$ we have
\begin{equation}
\sup_{\lambda,\gamma\in\R^d}\Big\vert\sqrt{Q_M (\lambda
+\gamma)}\sqrt{Q_M(\gamma)}-\sqrt{Q (\lambda +\gamma)}
\sqrt{Q(\gamma)}\Big\vert <\epsilon.\label{ap.23}
\end{equation} for sufficiently large $M$. Therefore,
\begin{eqnarray} &&
\bigg\{\int_{\{\vert \lambda\vert\le a\}}d\lambda
\bigg[\int_{\R^d}\sqrt{Q_M (\lambda +\gamma)}
\sqrt{Q_M(\gamma)}g(\lambda +\gamma)
g(\gamma)d\gamma\bigg]^p\bigg\}^{1/p}\label{ap.24}\\ &&
\le \epsilon\bigg\{\int_{\{\vert \lambda\vert\le a\}}d\lambda
\bigg[\int_{\R^d}g(\lambda +\gamma)
g(\gamma)d\gamma\bigg]^p\bigg\}^{1/p} \nonumber
\\ && +\bigg\{\int_{\{\vert \lambda\vert\le a\}}d\lambda
\bigg[\int_{\R^d}\sqrt{Q (\lambda +\gamma)}
\sqrt{Q(\gamma)}g(\lambda +\gamma)
g(\gamma)d\gamma\bigg]^p\bigg\}^{1/p}.\nonumber
\end{eqnarray} 
Also, since $\|g\|_{2}=1$
\begin{equation}
\int_{\{\vert \lambda\vert\le a\}}d\lambda
\bigg[\int_{\R^d}g(\lambda +\gamma) g(\gamma)d\gamma\bigg]^p\le
C_da^d\label{ap.25}
\end{equation} where $C_d$ is the volume of a $d$-dimensional unit ball. (\ref{ap.22}) then follows using the uniform continuity of $\wp_{\al,0}\ast \wp_{\al, 0}(\lambda )$, and finally  (\ref{64.2a}).
\qed

\section{A variational formula}\label{sec-varia}

The goal of this section is to prove the Theorem \ref{theo-varia}.
We begin with the following Lemma.

\begin{lemma}\label{lem-soboldir} Let $p\bb>d-\si>0$.
 Then 
\begin{equation}
\La_{d-\si}=:\sup_{g\in{\cal F}_\beta}\bigg\{ 
\bigg(\int_{ (R^{ d})^{ p}} {\prod_{ j=1}^{ p}g^{ 2}( x_{ j})\over |x_{ 1}+\cdots+x_{
p}|^{ d-\si}}
\prod_{ j=1}^{ p}\,d x_{ j}\bigg)^{1/p}-{\cal E}_\bb (g,g)\bigg\}<\ff.\label{67.0}
\end{equation}
\end{lemma}

Proof of Lemma \ref{lem-soboldir}:  By (\ref{60.2})
\bea &&
\bigg(\int_{ (R^{ d})^{ p}} {\prod_{ j=1}^{ p}g^{ 2}( x_{ j})\over |x_{ 1}+\cdots+x_{
p}|^{ d-\si}}
\prod_{ j=1}^{ p}\,d x_{ j}\bigg)^{1/p}\label{67.0j}\\
&&\hspace{ 1in}\leq
C\|g^{ 2}\|_{ pd/( ( p-1)d+\si)}=C\|g\|^{2}_{ 2pd/( ( p-1)d+\si)}.\nn
\eea

We then use the fact that for some $c<\ff$
\begin{equation}
\|f\|_{ 2pd/( ( p-1)d+\si)}\leq c\|\wh{f}\|_{ 2pd/(( p+1)d-\si)},\hspace{ .2in}f\in
\mathcal{S}( R^{ d})\label{67.0d}
\end{equation} and for any $r>0$
\begin{eqnarray} &&
\|\wh{f}\|_{ 2pd/(( p+1)d-\si)}^{ 2pd/(( p+1)d-\si)}\label{67.0e}\\
&& =\int_{ R^d}
{(r+|\la|^{ \bb})^{ pd/(( p+1)d-\si)} \over  (r+|\la|^{ \bb})^{pd/((
p+1)d-\si)}}|\wh{f}( \la)|^{2pd/(( p+1)d-\si)}\,d\la
\nn\\ && 
\leq \|  (r+|\la|^{ \bb})^{ - pd/(( p+1)d-\si)}\|_{(( p+1)d-\si)/(d-\si)}\nn\\ &&
\hspace{ .5in}\|  (r+|\la|^{ \bb})^{ pd/(( p+1)d-\si)}   |\wh{f}(
\la)|^{2pd/(( p+1)d-\si)}\|_{(( p+1)d-\si)/pd}.
\nonumber
\end{eqnarray} Now if $\|f\|_{ 2}=1$ then
\bea &&
\|  (r+|\la|^{ \bb})^{ pd/(( p+1)d-\si)}   |\wh{f}(
\la)|^{2pd/(( p+1)d-\si)}\|_{(( p+1)d-\si)/pd}\label{67.0f}\\
&&
=\(r+\mathcal{E}_{\bb}( f,f)\)^{pd/(( p+1)d-\si) }\nn
\eea and 
\begin{eqnarray} && h_{p, r}=:\|  (r+|\la|^{ \bb})^{ - pd/(( p+1)d-\si)}\|_{((
p+1)d-\si)/(d-\si)}\nn\\ && =\(\int_{ R^d} {1 \over (r+|\la|^{
\bb})^{ pd/(d-\si)}}\,d\la\)^{(d-\si)/(( p+1)d-\si) }.
\label{67.0g}
\end{eqnarray} 
Since $p\bb>d-\si$ this is finite and $\lim_{r\rar\ff }h_{p, r}=0$. 
Together we have shown that
\begin{equation}\qquad
\bigg(\int_{ (R^{ d})^{ p}} {\prod_{ j=1}^{ p}g^{ 2}( x_{ j})\over |x_{ 1}+\cdots+x_{
p}|^{ d-\si}}
\prod_{ j=1}^{ p}\,d x_{ j}\bigg)^{1/p}\leq c
h^{(( p+1)d-\si)/pd}_{p,r}
\(r+\mathcal{E}_{\bb}( g,g)\). \label{67.0k}
\end{equation} Our Lemma follows on taking $r$ sufficiently large so that $c
h^{(( p+1)d-\si)/pd}_{p,r}\leq 1$.\qed

Let $\mathcal{H}$ be a Hilbert space with norm $\|f\|$. We say that a  (possibly
unbounded) functional $L$ on $\mathcal{H}$ is positively homogeneous of order $k$ if for any
$\la\in R^{ 1}$
 and $f\in\mathcal{H}$ 
\begin{equation}
L( \la f)=|\la|^{ k}L( f).\label{b6.1}
\end{equation}
The following simple Lemma will be very useful.
\begin{lemma}\label{lem-geneq}
Let $L,\wt{L}$ be positive and positively homogeneous functionals on $\mathcal{H}$ of order $2$. 
For any $\th>0$ let
\begin{eqnarray} 
\La(\th)&=&\sup_{ \|f\|=1}\(\th L( f)-\wt{L}( f)\)\label{b6.2}\\ &=&
\sup_{f\in\mathcal{H}} { \(\th L( f)-\wt{L}( f)\)\over \|f\|^{ 2}}\nonumber
\end{eqnarray}
and assume that $\La(\th)$ is continuous.
Let 
\begin{eqnarray} 
J&=&\sup_{ \|f\|^{ 2}+\wt{L}( f)=1}L( f)\label{b6.3}\\ &=&
\sup_{f\in\mathcal{H}} { L( f)\over \|f\|^{ 2}+\wt{L}( f)}\nonumber
\end{eqnarray}
and assume that $J<\ff$.
Then
\begin{equation}
\La\({ 1\over J}\)=1.\label{b6.4}
\end{equation}
\end{lemma}

{\bf  Proof of Lemma \ref{lem-geneq}: } Fix $\ep>0$ and choose $g\in\mathcal{H}$ wth 
$\|g\|^{ 2}+\wt{L}( g)=1$ such that
\begin{equation}
L( g)\geq J-\ep.\label{b6.5}
\end{equation}
Then
\begin{eqnarray} 
\La\({ 1\over (J-\ep)}\)&\geq &{ \((J-\ep)^{ -1} L( g)-\wt{L}( g)\)\over \|g\|^{ 2}}\label{b6.6}\\
&\geq & { \((J-\ep)^{ -1}(
J-\ep)-\wt{L}( g)\)\over 1-\wt{L}( g)}=1.\nn
\end{eqnarray} 
By the continuity of $\La(\th)$, on taking $\ep\rar 0$ we see that $\La\({ 1\over J}\)\geq 1.$

On the other hand, by (\ref{b6.3}), for any $f\in\mathcal{H}$  
\begin{equation}
L( f)\leq J \(\|f\|^{ 2}+\wt{L}( f)\)\label{b6.7}
\end{equation} 
so that
\begin{eqnarray} \La\({ 1\over J}\) &=&\sup_{ \|f\|=1}\(J^{ -1} L( f)-\wt{L}( f)\)
\label{b6.7a}\\ &\leq &\sup_{ \|f\|=1}\(J^{ -1} J \(\|f\|^{ 2}+\wt{L}( f)\)-\wt{L}( f)\)=1.  \nonumber
\end{eqnarray}\qed

{\bf  Proof of Theorem \ref{theo-varia}}

We take $\mathcal{H}=L^{ 2}( R^{ d},\,dx)$,  $\wt{L}( f)= \mathcal{E}_{ \bb}( f,f)=( 2\pi)^{ -d}\int |\wh{f}( \la)|^{ 2}\psi ( \la)\,d\la$ and
\begin{equation}
L( f)=\(\int_{\R^{ pd}}{\prod_{ j=1}^{ p}|f(x_{ j})|^{ 2}\over
|x_{ 1}+\cdots+x_{ p}|^{ \si}}\prod_{ j=1}^{ p}\,dx_{ j}\)^{ 1/p}.\label{b6.8}
\end{equation}
If $f_{ \ep}( x)=\ep^{ d/2}f(\ep x)$ then $\|f_{ \ep}\|_{ 2}=\|f\|_{ 2}$, 
\begin{equation}
L( f_{ \ep})=\(\ep^{ dp}\int_{\R^{ pd}}{\prod_{ j=1}^{ p}|f(\ep x_{ j})|^{ 2}\over
|x_{ 1}+\cdots+x_{ p}|^{ \si}}\prod_{ j=1}^{ p}\,dx_{ j}\)^{ 1/p}=\ep^{ \si/p}L( f).\label{b6.9}
\end{equation}
Furthermore, $\wh{(f_{ \ep})}( \la)=\ep^{- d/2}\wh{f}( \la/\ep)$ so that
\begin{equation}
\wt{L}( f_{ \ep})=\ep^{- d}( 2\pi)^{ -d}\int |\wh{f}( \la/\ep)|^{ 2}\psi ( \la)\,d\la=\ep^{ \bb}\wt{L}(
f).\label{b6.10}
\end{equation}
Thus
\begin{eqnarray} 
\La(\th)&=&\sup_{ \|f\|_{2}=1}\(\th L( f)-\wt{L}( f)\)\label{b6.11}\\ &=&
\sup_{ \|f\|_{2}=1}\(\th L( f_{ \ep})-\wt{L}( f_{ \ep})\)\nonumber\\ &=&
\sup_{ \|f\|_{2}=1}\(\th \ep^{ \si/p} L( f)-\ep^{ \bb}\wt{L}( f)\).\nonumber
\end{eqnarray}
Taking $\ep=\th^{ 1/( \bb-\si/p)}$ we see that
\begin{equation}
\La(\th)=\th^{ \bb/( \bb-\si/p)}\La(1)\label{b6.12}
\end{equation}
which shows that $\La(\th)$ is continuous and that we can write (\ref{b6.4}) as 
\begin{equation}
J=( \La(1))^{1- \si/p\bb}.\label{b6.1}
\end{equation}

Recall that
\begin{eqnarray} &&
\rho\label{b6.9z}=\sup_{\|f\|_{ 2}=1}\int_{\R^d}\bigg[\int_{\R^d}{f(\lambda +\gamma)f(\gamma)\over\sqrt{ 1+\psi(\lambda
+\gamma)}\sqrt{ 1+\psi(\gamma)}}d\gamma\bigg]^p\phi_{ d-\si}(\lambda)d\lambda.
\end{eqnarray}
Setting 
$f=g/\sqrt{1+\psi}$ and using the notation $Q=(1+\psi)^{ -1}$ we have that
\begin{eqnarray} 
\rho&=&\sup_{( g,Qg)=1}\int_{\R^d}\bigg[\int_{\R^d}(Qg)(\lambda
+\gamma)(Qg)(\gamma)d\gamma\bigg]^p\phi_{ d-\si}(\lambda)d\lambda\label{b6.9y}\\&=&
\sup_{( g,Qg)=1}\int_{\R^d}\bigg[(Qg)\ast (\widetilde{Qg})\bigg]^p(-\lambda)\phi_{
d-\si}(\lambda)d\lambda\nn
\end{eqnarray}
where $\widetilde{f}( \ga)=f(-\ga )$.
Then, using $\mathcal{F}$ to denote the Fourier transform on $R^{ d}$, by Parseval's identity,
which can be justified as in the proof of Lemma \ref{lem-exprep},
\begin{eqnarray} 
\rho&=&
\sup_{( g,Qg)=1}( 2\pi)^{ -d}\int_{\R^d}\mathcal{F}\bigg[(Qg)\ast
(\widetilde{Qg})\bigg]^p(x)\mathcal{F}\phi_{ d-\si}(x)dx.\label{b6.9x}
\end{eqnarray}
Using the facts that $\mathcal{F}(f\ast g)=\mathcal{F}(f)\mathcal{F}(g)$, $\mathcal{F}(f
g)=( 2\pi)^{ -d}\mathcal{F}(f)\ast\mathcal{F}(g)$ and (\ref{a1.4}), and using the notation
$f^{ \ast p}$ for the
$p$-fold convolution product of $f $ with itself we see that
\begin{eqnarray} 
\rho&=&
\sup_{( g,Qg)=1}( 2\pi)^{ -d( p+1)}\int_{\R^d}\bigg[|\wh{Qg}|^{ 2}\bigg]^{ \ast p}(x){1
\over |x|^{
\si}}dx\label{b6.9u}\\&=&
\sup_{( 2\pi)^{ d}\|h\|_{ 2}^{ 2}+( 2\pi)^{ d}\wt{L}( h)=1}( 2\pi)^{ d( p-1)}\int_{\R^d}|h^{ 2}|^{ \ast p}(x){1
\over |x|^{
\si}}dx\nn
\end{eqnarray}
where in the last line we set $h=( 2\pi)^{ -d}\wh{Qg}$ so that $\wt{g}=Q^{ -1}\wh{h}$
and therefore 
$( g,Qg)=(\wh{h},Q^{ -1}\wh{h})=(\wh{h},( 1+\psi)\wh{h})=( 2\pi)^{ d}\|h\|_{ 2}^{ 2}+( 2\pi)^{ d}\wt{L}( h).$
By a change of variables we see that
\begin{equation}
\rho=\sup_{\|h\|_{ 2}^{ 2}+\wt{L}( h)=1}( 2\pi)^{ -d}\int_{\R^{ pd}}{\prod_{ j=1}^{
p}|h(x_{ j})|^{ 2}\over |x_{ 1}+\cdots+x_{ p}|^{ \si}}\prod_{ j=1}^{ p}\,dx_{ j}=( 2\pi)^{ -d}J^{ p}\label{b6.9v}
\end{equation}
and consequently by (\ref{b6.1})
\begin{equation}
\rho=( 2\pi)^{ -d}( \La(1))^{p- \si/\bb}.\label{b6.9w}
\end{equation}
\qed

 \section{Large deviations for $\zeta^{\ast}\big([0,1]^p\big)$}\label{sec-sup}

By Theorem \ref{theo-1}, the non-trivial part of Theorem \ref{theo-sup} is the upper bound.

\begin{lemma}\label{lem-exphold}
For any  $( t_1,\cdots,t_p)$, $M<\ff$  and any $\ga>0$ sufficiently small so that $\si'=\si+\ga$  satisfies (\ref{a1.2}),
there is a $c=c(M,\de)>0$ such that
\begin{equation}
\sup_{x}\E\exp\bigg\{c\sup_{\scriptstyle y\in B(x,M)\atop\scriptstyle
y\not =x}\bigg({\big\vert(\zeta^y-\zeta^x)\big([0,t_1]
\times\cdots\times [0,t_p]\big)\big\vert\over 
\vert y-x\vert^{\ga}}\bigg)^{1/p}\bigg\}<\infty.\label{9.3f}
\end{equation}
and
\begin{equation}
\sup_{x}\E\exp\bigg\{c\sup_{\scriptstyle y\in B(x,M)\atop\scriptstyle
y\not =x}\bigg({\big\vert(\zeta^y-\zeta^x)\big([0,\tau_1]
\times\cdots\times [0,\tau_p]\big)\big\vert\over 
\vert y-x\vert^{\ga}}\bigg)^{1/p}\bigg\}<\infty.\label{9.3}
\end{equation}
\end{lemma}

{\bf  Proof of Lemma  \ref{lem-exphold}}
By (\ref{jc.8})
 there is a $C_0=C_0(\zeta,\psi,p)>0$ such that
\begin{equation}
\sup_{
y\not =z}\E\bigg\vert{(\zeta^y-\zeta^z)
\big([0,t_1]
\times\cdots\times [0,t_p]\big)\over \vert y-z\vert^{\ga}}\bigg\vert^n
\le (n!)^p C_0^n\hskip.2in n=0,1,2,\cdots.
\label{8.32}
\end{equation}  

Recall that a function $\Psi$: $\R^+\longrightarrow \R^+$ is called a Young's function if
it is convex, increasing and satisfies $\Psi(0)=0$, 
$\displaystyle\lim_{x\to\infty}\Psi(x)=\infty$. The
Orlicz space ${\cal L}_\Psi (\Omega, {\cal A}, \P)$ is defined
as the linear space of all random  variables $X$ on the probability space
$(\Omega, {\cal A}, \P)$ such that
\begin{equation}
\vert\vert X\vert\vert_\Psi =\inf\big\{c>0;\hskip.1in 
\E\Psi (c^{-1}\vert X\vert)\le 1\big\}.
\label{8.33}
\end{equation}
It is known that $\vert\vert \cdot\vert\vert_\Psi$ defines a
norm (called the Orlicz norm) and ${\cal L}_\Psi (\Omega, {\cal A}, \P)$
becomes a Banach space under $\vert\vert \cdot\vert\vert_\Psi$.

We now choose the Young function $\Psi$  such that 
$\Psi(x)\sim\exp\big\{x^{1/p}\big\}$ as $x\to\infty$. By (\ref{8.32})
there is $c=c(\zeta, d,p)>0$ such that
\begin{equation}
\vert\vert (\zeta^y-\zeta^z)
\big([0,t_1]
\times\cdots\times [0,t_p]\big)
\vert\vert_\Psi\le c \vert y-z\vert^{\ga},\hspace{.2 in}\forall y,z.\label{8.34}
\end{equation}
By a standard chaining argument (see, e.g., \cite[Lemma 9]{CLR}), for any $\ga'<\ga$, $M<\ff$, uniformly in $x$
\begin{equation}
\bigg\vert\bigg\vert\sup_{\scriptstyle y\in B(x,M)\atop\scriptstyle
y\not =x}{\big\vert(\zeta^y-\zeta^x)\big([0,t_1]
\times\cdots\times [0,t_p]\big)\big\vert\over 
\vert y-x\vert^{\ga'}}\bigg\vert
\bigg\vert_\Psi<\infty\label{8.35}
\end{equation}
which leads to (\ref{9.3f}), after renaming $\ga'$ as $\ga$. The proof of (\ref{9.3}) is similar, as one can easily see that  (\ref{jc.8}) holds with all $t_{i}$ replaced by $\tau_{i}$.\qed

Now choosing $\ga$ so that (\ref{9.3f}) holds, pick $\la$ so that $(1+\ga\la)/p=\bb/\si$. By (\ref{a1.2}) we have that $\la>0$. It then follows from (\ref{9.3f}) that for some $C<\ff$ and all $t\geq 1$
\begin{eqnarray}
&&   \sup_{x\in\R^d} P\Big\{\sup_{y\in B(x,\epsilon t^{-\la})}
 \big\vert \zeta^x([0,1]^p)-\zeta^y([0,1]^p)\big\vert\ge \delta t\Big\}
\label{mod.10}\\
&&  \leq  \sup_{x\in\R^d} P\Big\{\sup_{y\in B(x,\epsilon t^{-\la})}
  { \big\vert \zeta^x([0,1]^p)-\zeta^y([0,1]^p)\big\vert  \over |x-y|^{\ga}}\ge { \delta t\over \epsilon^{\ga} t^{-\ga\la}}\Big\} \nonumber\\
&&  \leq Ce^{-\({ \delta t\over \epsilon^{\ga} t^{-\ga\la}}\)^{1/p}}= Ce^{-\({\de \over \ep^{\ga}}\)^{1/p}
t^{(1+\ga\la)/p}}= Ce^{-\({\de \over \ep^{\ga}}\)^{1/p}
t^{\bb/\si}}. \nonumber
\end{eqnarray}
Consequently,
\begin{eqnarray}
&&
\lim_{\epsilon\to 0^+}\limsup_{t\to\infty}t^{-\bb/\sigma}\log \sup_{x\in\R^d}\label{mod.11}\\
&&\hspace{.5 in}
P\Big\{\sup_{y\in B(x,\epsilon t^{-\la})}
 \big\vert \zeta^x([0,1]^p)-\zeta^y([0,1]^p)\big\vert\ge \delta t\Big\}
=-\infty.
 \nonumber
\end{eqnarray}

We first consider the case of $\bb=2$, the case of Brownian motion. 
By (\ref{mod.11}), for some $\la>0$ we have that  for 
any $\delta>0$ 
\begin{eqnarray}
&&
\lim_{\epsilon\to 0^+}\limsup_{t\to\infty}t^{-2/\sigma}\log \sup_{x\in\R^d}\label{xc.2}\\
&&\hspace{.5 in}
P\Big\{\sup_{y\in B(x,\epsilon t^{-\la})}
 \big\vert \zeta^x([0,1]^p)-\zeta^y([0,1]^p)\big\vert\ge \delta t\Big\}
=-\infty.
 \nonumber
\end{eqnarray}

Since the supremum of a Gaussian process has Gaussian tails, we have
\begin{equation}
\lim_{M\to\infty}\limsup_{t\to\infty}t^{-2/\sigma}\log
P\Big\{\sup_{s_1,\cdots, s_p\le 1}\vert X_1(s_1)+\cdots +
X_p(s_p)\vert\ge Mt^{1/\sigma}\Big\}=-\infty.
\label{xc.3}
\end{equation}

When $\sup_{s_1,\cdots, s_p\le 1}\vert X_1(s_1)+\cdots +
X_p(s_p)\vert\le Mt^{1/\sigma}$ and $\vert x\vert\ge 2Mt^{1/\sigma}$,
the quantity 
\begin{equation}
\zeta^x([0,1]^p)\le \int_0^1\!\!\cdots\!\!\int_0^1\Big(\vert x\vert -
\vert X_1(s_1)+\cdots +X_p(s_p)\vert\Big)^{-\sigma}ds_1\cdots ds_p
\label{xc.4}
\end{equation}
is bounded (so it is less than $t$). Thus
\begin{eqnarray}
&&
  P\Big\{
 \sup_{x\in\R^d}\zeta^x([0,1]^p)\ge t\Big\}\label{xc.5}\\
&&
\le P\Big\{
 \sup_{\vert x\vert\le 2Mt^{1/\sigma}}\zeta^x([0,1]^p)\ge t\Big\}\nn\\
&&
+P\Big\{\sup_{s_1,\cdots, s_p\le 1}\vert X_1(s_1)+\cdots +
X_p(s_p)\vert\ge Mt^{1/\sigma}\Big\}.
 \nonumber
\end{eqnarray}

The cardinality of an $\epsilon t^{-\la}$-net on the ball of radius $2Mt^{1/\sigma}$ 
is of the order $O(t^{d(\la+\si^{-1})})$. This gives  
\begin{eqnarray}
&&
P\Big\{
 \sup_{\vert x\vert\le 2Mt^{1/\sigma}}\zeta^x([0,1]^p)\ge t\Big\}\label{xc.6}\\
&&
\le Ct^{d(\la+\si^{-1})}\lc
\sup_{x\in\R^d}P\Big\{\zeta^x([0,1]^p)\ge (1-\delta)t\Big\}\right.\nn\\
&&\left.
+\sup_{x\in\R^d} P\Big\{\sup_{y\in B(x,\epsilon t^{-\la})}
 \big\vert \zeta^x([0,1]^p)-\zeta^y([0,1]^p)\big\vert\ge \delta t\Big\}\rc.\nn
\end{eqnarray}

Summarizing what we have
\begin{eqnarray}
&&
\limsup_{t\to\infty}t^{-2/\sigma}\log P\Big\{
\sup_{x\in\R^d} \zeta^x([0,1]^p)\ge t\Big\}\label{xc.7}\\
&&\le \max\bigg\{\limsup_{t\to\infty}t^{-2/\sigma}\log\sup_{x\in\R^d}P\Big\{
\zeta^x([0,1]^p)\ge (1-\delta)t\Big\},\nn\\
&&
\limsup_{t\to\infty}t^{-2/\sigma}\log \sup_{x\in\R^d}
P\Big\{\sup_{y\in B(x,\epsilon t^{-\la})}
 \big\vert \zeta^x([0,1]^p)-\zeta^y([0,1]^p)\big\vert\ge \delta t\Big\},\nn\\
&&
\limsup_{t\to\infty}t^{-2/\sigma}\log
P\Big\{\sup_{s_1,\cdots, s_p\le 1}\vert X_1(s_1)+\cdots +
X_p(s_p)\vert\ge Mt^{1/\sigma}\Big\}\bigg\}.\nn
\end{eqnarray}
Letting $M\to\infty$ and $\epsilon\to 0$ we have
\begin{eqnarray}
&&
\limsup_{t\to\infty}t^{-2/\sigma}\log P\Big\{
\sup_{x\in\R^d} \zeta^x([0,1]^p)\ge t\Big\}\label{xc.8}\\
&&
\le\limsup_{t\to\infty}t^{-2/\sigma}\log\sup_{x\in\R^d}\P\Big\{
\zeta^x([0,1]^p)\ge (1-\delta)t\Big\}.\nn
\end{eqnarray}

Using (\ref{domin}) and (\ref{a2.15a})
\begin{equation}
\limsup_{n\to\infty}{1\over n}\log { 1\over (n!)^{\si/\bb}}\sup_{x\in\R^d}
\E\Big[\zeta^{x} ([0,1]^p)^n\Big]\le
\log \Big({\bb p\over \bb p -\si}\Big)^{\bb p -\si\over\bb} +\log \rho.\label{a2.15aj}
\end{equation}
The (easy part of the) proof of \cite[Lemma 2.3]{KM} then shows that
\begin{eqnarray}
&&
\limsup_{t\to\infty}t^{-\bb/\sigma}\log\sup_{x\in\R^d}P\Big\{
\zeta^x([0,1]^p)\ge  t\Big\}\le  -
{\sigma\over \bb} \Big({p\bb  -\sigma\over p\bb}
\Big)^{p\bb -\sigma\over\sigma}\rho^{- \bb/\sigma}.
\label{xc.10}
\end{eqnarray}
With $\bb=2$ we have
\begin{eqnarray}
&&
\limsup_{t\to\infty}t^{-2/\sigma}\log\sup_{x\in\R^d}P\Big\{
\zeta^x([0,1]^p)\ge (1-\delta)t\Big\}\label{xc.10j}\\
&&\hspace{1 in}\le  -
(1-\delta)^{2/\sigma} 
{\sigma\over 2} \Big({2p -\sigma\over 2p}
\Big)^{2p -\sigma\over\sigma}\rho^{-2/\sigma}.
\nn
\end{eqnarray}
Thus
\begin{eqnarray}
&&
\limsup_{t\to\infty}t^{-2/\sigma}\log P\Big\{
\sup_{x\in\R^d} \zeta^x([0,1]^p)\ge t\Big\}\label{xc.11}\\
&&\hspace{1 in}
\le  -
(1-\delta)^{2/\sigma} 
{\sigma\over 2} \Big({2p -\sigma\over 2p}
\Big)^{2p -\sigma\over\sigma}\rho^{-2/\sigma}.
\nn
\end{eqnarray}
Letting $\delta\to 0^+$ gives (\ref{a1.28sup}).

We now consider $\bb\neq 2$. 
We will show  that there exists $c_1>0$ such that
\begin{equation}
E\Big(\exp\( c_1 \lc\sup_{z\in R^d}\zeta^x\big([0,\tau_1]
\times\cdots\times [0,\tau_p]\big)\rc^{1/p}\) \Big)<\infty.\label{RB-E1}
\end{equation}
It will follow from this that for some $c_{2}<\ff$
\begin{equation}
E\Big( \lc\sup_{z\in R^d}\zeta^x\big([0,\tau_1]
\times\cdots\times [0,\tau_p]\big)\rc^{n/p} \Big)\leq n!c^{n}_{2}\label{RB-E1a}
\end{equation}
for all $n$. Hence, taking $n=mp$
\begin{equation}
E\Big( \lc\sup_{z\in R^d}\zeta^x\big([0,\tau_1]
\times\cdots\times [0,\tau_p]\big)\rc^{m} \Big)\leq (mp)!c^{pm}_{2}
\leq (m!)^{p}c^{m}_{3}.\label{RB-E1b}
\end{equation}
Using (\ref{a2.15}) and Stirling's formula as in (\ref{a2.15a}) we obtain 
\begin{equation}
\limsup_{n\to\infty}{1\over n}\log { 1\over (n!)^{\si/\bb}}
\E\Big[\zeta^{\ast} ([0,1]^p)^n\Big]\le c_{4}<\ff.\label{a2.15arx}
\end{equation}
Then once again the  (easy part of the) proof of \cite[Lemma 2.3]{KM} will show  
that for some $0<C<\ff$
\begin{eqnarray}
&&
\limsup_{t\to\infty}t^{-\bb/\sigma}\log P\Big\{
\zeta^{\ast} ([0,1]^p)\ge  t\Big\}\le  -C.
\label{xc.10w}
\end{eqnarray}
Thus it only remains to show (\ref{RB-E1}).

\begin{lemma}\label{RB-L1}
Let $X_t$ be a $d$-dimensional symmetric stable process
of order $\beta$ and $\tau$ an independent exponential of
parameter 1. Then   there exists a constant $c_1$ such that for $D>0$,
\begin{equation} P(\sup_{s\leq \tau}|X_s|\geq D)\leq \frac{c_1}
{D^\beta}.\label{RB-E2}
\end{equation}
\end{lemma}

\noindent {\bf Proof.} It is well known, \cite[Proposition 2.2]{Kolokoltsov}, that the density of
$X_t$ satisfies
$$p(t,x,y)\leq ct/|x-y|^{d+\beta}.$$
(A better estimate is possible for larger $t$, but this is not
needed.)  Integrating over $|y-x|\geq D$, we obtain
\begin{equation}\label{RB-E3}
P^x(|X_t-X_0|\geq D)\leq \frac{ct}{D^\beta}.
\end{equation}

We now obtain an estimate
on the exit probabilities. Let $S=\inf\{s: |X_s|\geq D\}$. If
$\sup_{s\leq t} |X_s|\geq D$, then $S\leq t$ and either $|X_t|\geq D/2$ or
$|X_t|\leq D/2$, so that $|X_S-X_t|\geq D/2$. Thus  
\begin{align*}
P(\sup_{s\leq t} |X_s|\geq D)& 
\leq P(|X_t|\geq D/2)\\
&~~~~+P(S<t, |X_t-X_S|\geq D/2).
\end{align*}
The first term on the right is bounded by 
$ct/D^\beta$ using (\ref{RB-E3}). The second term on the
right is bounded by
$$\int_0^t P(|X_t-X_s|\geq D/2)\, P(S\in ds)\leq 2c\int_0^t (t-s)/D^\beta
P(S\in ds) \leq 2ct/D^\beta$$
using (\ref{RB-E3}) again and the Markov property of $X$.

Finally,
\begin{align*}
P(\sup_{s\leq \tau} |X_s|\geq D)&=\int_0^\infty e^{-t} P(\sup_{s\leq t}
|X_s|\geq D) \, dt\leq \int_0^\infty e^{-t}\frac{ct}{D^\beta}\, dt\\
&\leq c/D^\beta
\end{align*}
as desired.
\qed

\begin{lemma}\label{RB-L2} 
Suppose for each $z\in R^d$ there is a random variable $Y^z$
such that $z\to Y^z$ is continuous, a.s., and there exist $\delta$, $A$ and $B$
such that
\begin{align}
E e^{A|Y^z|}&\leq B, \qquad z\in R^d, \label{RB-E5a}\\
Ee^{A|Y^z-Y^{z'}|/|z-z'|^\delta}&\leq B, \qquad
z,z'\in R^d.\label{RB-E5b}
\end{align}
Then there exist $c_1$ and $c_2$ such that
for every $M\geq 1$
\begin{equation}\label{RB-E4}
E \exp\Big(c_1A \sup_{|z|\leq M} |Y^z|\Big)\leq c_2M^{2d} B.
\end{equation}
\end{lemma}

\noindent {\bf Proof.} Let $Q_k=B(0,M)\cap 2^{-k}Z^d$ and $Q=\cup_k Q_k$.
Since $z\to Y^z$ is continuous, it suffices to bound
\begin{equation}\label{RB-E6}
E \exp\Big(c_1A \sup_{|z|\in Q} |Y^z|\Big).
\end{equation}
If $z\in Q$, we write
$$z=z_0+(z_1-z_0)+(z_2-z_1)+\cdots.$$
Here $z_i$ is the point of $Q_i$ closest to $z$, with some convention 
for breaking ties. Since $z\in Q_k$ for some $k$, the above sum
is actually a finite one.

If $|Y^z|\geq \lam$, then either the event $R$ holds: $|Y^{z_0}|\geq \lam/2$ for some
$z_0\in Q_0$, or for some $i$ the event $S_{i}$ holds:  $|Y^{z_{i+1}}-Y^{z_i}|\geq \lam/20i^2$
for some pair $z_i, z_{i+1}$ with $z_i\in Q_i$, $z_{i+1}\in Q_{i+1}$,
and $|z_i-z_{i+1}| \leq \sqrt d 2^{-i}$.

Since there are at most $M^d$ points in $Q_0$,
using (\ref{RB-E5a}) we see the probability of the event $R$ 
is bounded by
$$cBM^d e^{-A\lam/2}.$$

For each $i$, there are at most $cM^{2d}2^{c'id}$ pairs $z_i, z_{i+1}$
as in the definition of the even $S_i$, so the probability of the event $S_i$
is bounded by
\begin{align}
cBM^{2d}2^{c'id}&\sup_{|z_i-z_{i+1}|\leq \sqrt d2^{-i}}
P\Big(\frac{|Y^{z_{i+1}}-Y^{z_i}|}{|z_{i+1}-z_i|^\delta}
\geq \frac{\lam A}{20 i^2 (\sqrt d 2^{-i})^\delta}\Big)\\
&\leq cBM^{2d} 2^{c'id} \exp\Big( -c''\frac{\lam A}{20i^2 (\sqrt d2^{-i})^\delta}
\Big).
\end{align}
If we sum over $i$ the probabilities of the events $S_i$  holding
and add to that the probability of the event $R$ holding,
we obtain
$$P(\sup_{z\in Q} |Y^z|\geq \lam)\leq cB M^{2d} e^{-c'\lam A}.$$
Our result follows from this.
\qed

We now prove (\ref{RB-E1}). Let $X^i_t$, $i=1, \ldots, p$,  be independent $d$-dimensional symmetric stable processes of order        
      $\beta$. We write simply $\zeta^z$ for
$\zeta^z([0, \tau_1]\times \cdots \times [0, \tau_p])$
and $Z_i$ for $\sup_{s\leq \tau_{i}} |X^i_s|$. 
We will choose $c_1$ later. 

It follows from (\ref{a2.8}) and (\ref{a2.10}) that there exists $c_2$ such that
\begin{equation}
\sup_{z\in R^d}E \exp\Big( c_2 |\zeta^z|^{1/p}\Big)<\infty,\label{RB-E1a}
\end{equation}
and using (\ref{9.3}) and the fact that $|a^{1/p}-b^{1/p}|\leq |a-b|^{1/p}$, we can choose $c_2$ such that also
\begin{equation}
\sup_{z,z'\in R^d} E \exp\Big(c_2 |\,|\zeta^z|^{1/p}-|\zeta^{z'}|^{1/p}| 
/|z-z'|^\delta\Big)<\infty.\label{RB-E1b}
\end{equation}

Write
\begin{align}
E e^{c_1\sup_z |\zeta^z|^{1/p}}
&=E\Big[e^{c_1\sup_z |\zeta^z|^{1/p}}; \max_{1\leq i\leq p} 
Z^i    
\leq 1\Big]\\
&~~~~+\sum_{k=0}^\infty E\Big[ e^{c_1\sup_z |\zeta^z|^{1/p}}; 2^k\leq
\max_{1\leq i\leq p} Z_i\leq 2^{k+1}\Big]\\
&:= I+\sum_{k=0}^\infty J_k.
\end{align}
Now write 
\begin{equation}\label{RB-E7}
I\leq E\Big[ e^{c_1 \sup_{|z|\leq 2p}|\zeta^z|^{1/p}}\Big]
+E\Big[ e^{c_1 \sup_{|z|> 2p}|\zeta^z|^{1/p}}; \max_{1\leq i\leq p} 
Z \leq 1\Big]=I'+I''.
\end{equation}
Provided $c_1<c_2$, then 
 $I'$ is finite by Lemma \ref{RB-L2} with $Y^z=|\zeta^z|^{1/p}$. 
If $|z|\geq 2p$ and $\max_{1\leq i\leq p} Z_i\leq 1$, then 
$$|Z_1+\cdots +Z_p-z|\geq p,$$
and hence 
$$I''\leq E e^{c_1p^{-\delta/p}(\tau_1\cdots \tau_p)^{1/p}}
=\int_{R^{p}_{+}}e^{-\sum_{j=1}^{p}t_{j}}e^{c_1p^{-\delta/p}(t_1\cdots t_p)^{1/p}}
\,dt_{1}\cdots\,dt_{p}.$$
Since $(t_1\cdots t_p)^{1/p}\leq \max_{1\leq j\leq p}t_j\leq \sum_{j=1}^{p}t_{j}$ we see that
\begin{equation}\label{RB-E9}
E e^{c(\tau_1\cdots \tau_p)^{1/p}}<\infty
\end{equation}
if $c$ is small enough.

Combining with the estimate for $I'$ shows that $I$ is finite, provided
$c_1<c_2$ and $c_2$ is sufficiently small.

We turn to $J_k$ and write
\begin{align*} J_k&\leq E\Big[e^{c_1\sup_{|z|\leq p2^{k+1}} |\zeta^z|^{1/p}}
; 2^k\leq \max_{1\leq i\leq p} Z_i\leq 2^{k+1}\Big]\\
&~~~~+ E\Big[e^{c_1\sup_{|z|> p2^{k+1}} |\zeta^z|^{1/p}}
; 2^k\leq \max_{1\leq i\leq p} Z_i\leq 2^{k+1}\Big]\\
&=J'_k+J_k''.
\end{align*}

For $J_k'$ we apply H\"older's inequality with $\frac1{r}+\frac1{s}$
and $r$ and $s$ to be chosen later.
Then  
\begin{align*}
J'_k&\leq \Big(E e^{c_1r \sup_{|z|\leq p2^{k+1}} |\zeta^z|^{1/p}}\Big)^{1/r}
\Big(P(\max_{1\leq i\leq p} Z_i\geq 2^{k})\Big)^{1/s}\\
&\leq \Big(c 2^{2dk}\Big)^{1/r} \Big(\frac{c}{2^{k\beta}}\Big)^{1/s}
\end{align*}
by Lemma \ref{RB-L2}. 
We now choose $r$ and $s$ so that $\eta:=\beta /s- 2d/r>0$, and hence
$2^{k\beta/s}\geq 2^{\eta k}  2^{2dk/r}$.
This proves $J_k'$ is summable in $k$.

To handle $J_k''$, if $\max_{1\leq i\leq p} Z_i\leq 2^{k+1}$ and
$|z|\geq p2^{k+2}$, then 
$$\sup_{s_1\leq 1, \ldots, s_p\leq 1}
|X^1_{s_1}+\cdots+X^p_{s_p}-z|\geq p2^{k+1},$$
and hence
$$J_k''\leq  E \Big[ e^{c_1(p2^{k+1})^{-\sigma/p}(\tau_1\cdots \tau_p)^{1/p}};
\max_{1\leq i\leq p} Z_i\geq 2^k\Big].$$
Using Cauchy-Schwarz, we obtain
$$J_k''\leq e^{c_1(p2^{k+1})^{-\sigma/p}}\,\,
P(\max_{1\leq i\leq p} Z_i\geq 2^k).
$$
The second factor is less than or equal to $c/2^{k\beta}$, which
is summable in $k$.

Finally we choose $c_1$ small enough that $c_1r<c_2$, and the proof  of (\ref{RB-E1})
 complete.\qed

\section{Laws of the iterated logarithm}\label{sec-lil}

The upper bound in (\ref{a1.28lils}) and therfore the upper bound in (\ref{a1.28lil})
follows from Theorem \ref{theo-1}, the  scaling property given in (\ref{a1.3}), and
a standard procedure using the Borel-Cantelli lemma. It remains to
prove that  
\begin{eqnarray}
&&\limsup_{t\to\infty}t^{-{p\bb-\sigma \over \bb}}(\log\log t)^{-\si/\bb}
\zeta ([0,t]^p)\geq \({\sigma\over\bb}\)^{-\si/\bb}
\Big({p\bb -\sigma\over p\bb}
\Big)^{\sigma-p\bb \over\bb}\rho \label{9.1} 
\end{eqnarray}
almost surely.

We first prove that
\begin{equation}
\lim_{\delta\to 0^+}\liminf_{t\to\infty} t^{-1}\log 
\P\Big\{\inf_{\vert y\vert\le\delta}\zeta^y\big([0,t]^p\big)\ge t^p\Big\}
\ge - {\sigma\over\bb}
\Big({p\bb -\sigma\over p\bb}
\Big)^{p\bb -\sigma\over\sigma}\rho^{-\bb/\sigma}.\label{9.2}
\end{equation}

Using (\ref{9.3}) and  Chebyshev's inequality we have that for any $\epsilon >0$,
\begin{equation}
\limsup_{\delta\to 0^+}\limsup_{t\to\infty} t^{-1}\log\P\Big\{
\sup_{\vert y\vert\le\delta}\big\vert(\zeta^0-\zeta^y)
\big([0,\tau_1]\times\cdots\times [0,\tau_p]\big)\big\vert
\ge \epsilon t^p\Big\}=-\infty.
\label{9.4}
\end{equation}

On the other hand,
\begin{eqnarray}
&&\quad\P\Big\{
\sup_{\vert y\vert\le\delta}\big\vert(\zeta^0-\zeta^y)
\big([0,\tau_1]\times\cdots\times [0,\tau_p]\big)\big\vert
\ge \epsilon t^p\Big\}\label{9.5}\\
&&=\int_0^\infty\!\!\cdots\!\!\int_0^\infty e^{-(t_1+\cdots +t_p)}\nonumber \\
&&\hspace{1 in}
\P\Big\{
\sup_{\vert y\vert\le\delta}\big\vert(\zeta^0-\zeta^y)
\big([0,t_1]\times\cdots\times [0,t_p]\big)\big\vert
\ge \epsilon t^p\Big\}dt_1\cdots dt_p \nonumber \\
&&\ge\int_{(1-\epsilon)t}^t\!\!\cdots\!\!\int_{(1-\epsilon)t}^t
e^{-(t_1+\cdots +t_p)}\nonumber \\
&&\hspace{1 in}
\P\Big\{
\sup_{\vert y\vert\le\delta}\big\vert(\zeta^0-\zeta^y)
\big([0,t_1]\times\cdots\times [0,t_p]\big)\big\vert
\ge \epsilon t^p\Big\}dt_1\cdots dt_p \nonumber \\
&&\ge \big(e^{-(1-\epsilon)t}-e^{-t}\big)^p\nonumber \\
&&\hspace{.5 in}
\inf_{(1-\epsilon)t\le t_1,\cdots, t_p\le t}
\P\Big\{
\sup_{\vert y\vert\le\delta}\big\vert(\zeta^0-\zeta^y)
\big([0,t_1]\times\cdots\times [0,t_p]\big)\big\vert
\ge \epsilon t^p\Big\}.\nonumber
\end{eqnarray}
So we have
\begin{eqnarray}
&&\limsup_{\delta\to 0^+}\limsup_{t\to\infty} t^{-1}\log
\inf_{(1-\epsilon)t\le t_1,\cdots, t_p\le t}\label{9.6}\\
&&\hskip.3in\P\Big\{
\sup_{\vert y\vert\le\delta}\big\vert(\zeta^0-\zeta^y)
\big([0,t_1]\times\cdots\times [0,t_p]\big)\big\vert
\ge \epsilon t^p\Big\}=-\infty.\nonumber
\end{eqnarray}

For any $t$ and $(1-\epsilon)t\le t_1,\cdots, t_p\le t$,
\begin{eqnarray}
&&\inf_{\vert y\vert\le\delta}\zeta^y\big([0,t]^p\big)
\ge \inf_{\vert y\vert\le\delta}\zeta^y\big([0,t_1]\times\cdots\times
[0,t_p]\big)\label{9.7}\\
&&\ge\zeta^0\big([0,t_1]\times\cdots\times
[0,t_p]\big)-\sup_{\vert y\vert\le\delta}\big\vert(\zeta^0-\zeta^y)
\big([0,t_1]\times\cdots\times
[0,t_p]\big)\big\vert \nonumber \\
&&\ge\zeta^0\big([0,(1-\epsilon)t]^p\big)
-\sup_{\vert y\vert\le\delta}\big\vert(\zeta^0-\zeta^y)
\big([0,t_1]\times\cdots\times
[0,t_p]\big)\big\vert.\nonumber
\end{eqnarray}
Hence,
\begin{eqnarray}
&&\P\Big\{\inf_{\vert x\vert\le\delta}\zeta^x\big([0,t]^p\big)
\ge t^p\Big\}\label{9.8}\\
&&+\inf_{(1-\epsilon)t\le t_1,\cdots, t_p\le t}
\P\Big\{
\sup_{\vert x\vert\le\delta}\big\vert(\zeta^0-\zeta^x)
\big([0,t_1]\times\cdots\times [0,t_p]\big)\big\vert
\ge \epsilon t^p\Big\} \nonumber \\
&&\ge \P\Big\{\zeta^0\big([0,(1-\epsilon)t]^p\big)\ge (1+\epsilon) t^p\Big\}.\nonumber
\end{eqnarray}
Consequently,
\begin{eqnarray}
&&\max\bigg\{\liminf_{t\to\infty} t^{-1}\log\P\Big\{
\inf_{\vert y\vert\le\delta}\zeta^y\big([0,t]^p\big)
\ge t^p\Big\},\label{9.9}\\
&&\hspace{.5 in}\limsup_{t\to\infty} t^{-1}\log\inf_{(1-\epsilon)t\le t_1,\cdots, t_p\le t}\nonumber \\
&&\hspace{1 in}
\P\Big\{
\sup_{\vert y\vert\le\delta}\big\vert(\zeta^0-\zeta^y)
\big([0,t_1]\times\cdots\times [0,t_p]\big)\big\vert
\geq \epsilon t^{p}\Big\}\bigg\} \nonumber \\
&&\ge\liminf_{t\to\infty} t^{-1}\log
\P\Big\{\zeta^0\big([0,(1-\epsilon)t]^p\big)\ge (1+\epsilon) t^p\Big\}.\nonumber
\end{eqnarray}
Notice that by the scaling (\ref{a1.3})
\begin{equation}
\P\Big\{\zeta^0\big([0,(1-\epsilon)t]^p\big)\ge (1+\epsilon) t^p\Big\}
=\P\Big\{\zeta^0\big([0,1]^p\big)\ge (1+\epsilon)
(1-\epsilon)^{-{\bb p-\si\over\bb}} 
t^{\si/\bb}\Big\},
\label{9.10}
\end{equation}
so that by  Theorem \ref{theo-1},
\begin{eqnarray}
&&\lim_{t\to\infty} t^{-1}\log
\P\Big\{\zeta^0\big([0,(1-\epsilon)t]^p\big)\ge (1+\epsilon) t^p\Big\}\label{9.11}\\
&&=-(1+\epsilon)^{\bb /\si}(1-\epsilon)^{-{p\bb -\sigma\over\sigma}}
{\sigma\over\bb}
\Big({p\bb -\sigma\over p\bb}
\Big)^{p\bb -\sigma\over\sigma}\rho^{-\bb/\sigma}.\nn
\end{eqnarray}
Let $\delta\to 0^+$ in (\ref{9.9}). By (\ref{9.9}), (\ref{9.6}) and (\ref{9.11}) we obtain
\begin{eqnarray}
&&\lim_{\delta\to 0^+}\liminf_{t\to\infty} t^{-1}\log\P\Big\{
\inf_{\vert y\vert\le\delta}\zeta^y\big([0,t]^p\big)
\ge t^p\Big\}\label{9.12}\\
&&\ge -(1+\epsilon)^{\bb /\si}(1-\epsilon)^{-{p\bb -\sigma\over\sigma}}
{\sigma\over\bb}
\Big({p\bb -\sigma\over p\bb}
\Big)^{p\bb -\sigma\over\sigma}\rho^{-\bb/\sigma}.\nonumber
\end{eqnarray}
Letting $\epsilon\to 0^+$ on the right hand side leads to (\ref{9.2}).
\medskip

We now come to the proof of (\ref{9.1}).
For each $k\ge 1$, write $t_k=k^k$ and define
\begin{equation}
X_{j,k}(t)=X_j(t_k+t)-X_j(t_k) \hskip.2in t\ge 0,\hskip.1in j=1,\cdots, p
,\hskip.1in k=1,2,\cdots.
\label{9.13}
\end{equation}
Let $\zeta_k^x([a,b]^{p})$ be the Riesz potential of the additive stable process
\begin{equation}
\ol{ X}_k(s_1,\cdots,s_p)=X_{1,k}(s_1)+\cdots +X_{p,k}(s_p).
\label{9.14}
\end{equation}
Then for each $k$,
$\{\zeta_k^{x}\,,x\in R^{d}\}\buildrel d\over =\{\zeta^{x}\,,x\in R^{d}\}$.
\medskip
Let $\delta >0$ be a small number which will be specified later.
Write $Y_k=X_1(t_k)+\cdots +X_p(t_k)$. A rough estimate gives that
with probability 1, the inequality
\begin{equation}
\vert Y_k\vert\le 2^{-1}\delta\Big({t_{k+1}\over\log\log t_{k+1}}
\Big)^{1/\bb}
\label{9.15}
\end{equation}
eventually holds. Therefore 
\begin{eqnarray}
&&\zeta \big([t_k,t_{k+1}]^p\big)
=\zeta_k^{ Y_k}\big([0,t_{k+1}-t_k]^p\big)\label{9.16}\\
&&\ge \inf_{\vert y\vert\le \delta(t_{k+1}/\log\log t_{k+1})^{1/\beta}}
\zeta_k^y\big([0,t_{k+1}-t_k]^p\big)\nonumber
\end{eqnarray}
eventually holds, almost surely.
\medskip
For each $k$, by the scaling (\ref{a1.3}),
\begin{eqnarray}
&&\inf_{\vert y\vert\le \delta(t_{k+1}/\log\log t_{k+1})^{1/\beta}}
\zeta_k^y\big([0,t_{k+1}-t_k]^p\big)\label{9.17}\\
&&\buildrel d\over =
\inf_{\vert y\vert\le \delta(t_{k+1}/\log\log t_{k+1})^{1/\beta}}
\zeta^y\big([0,t_{k+1}-t_k]^p\big) \nonumber \\
&&\buildrel d\over =
\Big({t_{k+1}\over \log\log t_{k+1}}\Big)^{\beta p-\si\over\beta}
\inf_{\vert y\vert\le\delta}\zeta^y\big([0,\hskip.05in
t_{k+1}^{-1}(t_{k+1}-t_k)\log\log t_{k+1}]^p\big).\nonumber
\end{eqnarray}
Let $\theta >0$ satisfy
\begin{equation}
\theta <\Big({\beta\over \si}\Big)^{\si\over\beta}
\Big(1-{\si\over\beta p}\Big)^{-(p-{\si\over\beta})}\rho.
\label{8.34}
\end{equation}
We have
\begin{eqnarray}
&&
\P\Big\{\inf_{\vert y\vert\le \delta(t_{k+1}/\log\log t_{k+1})^{1/\beta}}
\zeta_k^y\big([0,t_{k+1}-t_k]^p\big)\ge \theta
t_{k+1}^{\beta p-\si\over\beta}(\log\log t_{k+1})^{\si\over\beta}\Big\}\label{9.18}\\
&&=\P\Big\{\inf_{\vert y\vert\le\delta}\zeta^y\big([0,\hskip.05in
t_{k+1}^{-1}(t_{k+1}-t_k)\log\log t_{k+1}]^p\big)\ge 
\theta (\log\log t_{k+1})^p\Big\}.\nonumber
\end{eqnarray}
Using the scaling (\ref{a1.3}) once again 
\begin{eqnarray}
&&\qquad\P\Big\{\inf_{\vert y\vert\le\delta}\zeta^y\big([0,\hskip.05in
t_{k+1}^{-1}(t_{k+1}-t_k)\log\log t_{k+1}]^p\big)\ge 
\theta (\log\log t_{k+1})^p\Big\}
\label{9.18a}\\
&& =\P\Big\{\inf_{\vert y\vert\le\delta \theta^{1/\si}}\zeta^y\big([0,\hskip.05in
t_{k+1}^{-1}(t_{k+1}-t_k)\theta^{\bb/\si}\log\log t_{k+1}]^p\big)\ge 
 (\theta^{\bb/\si}\log\log t_{k+1})^p\Big\}.  \nonumber
\end{eqnarray}
By (\ref{9.2}), therefore, one can take $\delta >0$ sufficiently small so that
\begin{eqnarray}
&&\liminf_{k\to\infty}{1\over\log\log t_{k+1}}\log
\P\Big\{\inf_{\vert y\vert\le \delta(t_{k+1}/\log\log t_{k+1})^{1/\beta}}
\zeta_k^y\big([0,t_{k+1}-t_k]^p\big)\label{9.19}\\
&&\hskip.6in\ge \theta
t_{k+1}^{\beta p-\si\over\beta}(\log\log t_{k+1})^{\si\over\beta}\Big\}>-1.\nonumber
\end{eqnarray}
Consequently,
\begin{equation}
\sum_k\P\Big\{\inf_{\vert y\vert\le \delta(t_{k+1}/
\log\log t_{k+1})^{1/\beta}}
\zeta_k^y\big([0,t_{k+1}-t_k]^p\big)\ge \theta
t_{k+1}^{\beta p-\si\over\beta}(\log\log t_{k+1})^{\si\over\beta}\Big\}=\infty.
\label{9.20}
\end{equation}
Notice that
\begin{equation}
\inf_{\vert y\vert\le \delta(t_{k+1}/
\log\log t_{k+1})^{1/\beta}}
\zeta_k^y\big([0,t_{k+1}-t_k]^p\big)\hskip.2in k=1,2,\cdots
\label{9.21}
\end{equation}
is an independent sequence. By the Borel-Cantelli lemma,
\begin{equation}
\limsup_{k\to\infty}t_{k+1}^{-{\beta p- \si\over\beta}}
(\log\log t_{k+1})^{-{\si\over\beta}}\inf_{\vert y\vert\le \delta(t_{k+1}/
\log\log t_{k+1})^{1/\beta}}
\zeta_k^y\big([0,t_{k+1}-t_k]^p\big)\ge\theta\hskip.2in a.s.
\label{9.22}
\end{equation}
By (\ref{9.16}),
\begin{equation}
\limsup_{k\to\infty}t_{k+1}^{-{\beta p- \si\over\beta}}
(\log\log t_{k+1})^{-{\si\over\beta}}\zeta \big([t_k,t_{k+1}]^p\big)
\ge\theta\hskip.2in a.s.
\label{9.23}
\end{equation}
Consequently,
\begin{equation}
\limsup_{t\to\infty}t^{-{\beta p-\si\over\beta}}(\log\log t)^{-{\si\over\beta}}
\zeta \big([0,t]^p\big)
\ge\theta \hskip.2in a.s.
\label{9.24}
\end{equation}
Letting
\begin{equation}
\theta \uparrow \Big({\beta\over \si}\Big)^{\si\over\beta}
\Big(1-{\si\over\beta p}\Big)^{-(p-{\si\over\beta})}\rho
\label{9.25}
\end{equation}
proves (\ref{9.1}). \qed

\section{Appendix: Sobolev-type inequalities}\label{sec-app}

\begin{lemma}\label{lem-multisob}For any $q>1$ and integer $p\geq 1$
\begin{equation}
\|f_{ 1}\ast\cdots\ast f_{p}\|_{ q}\leq C^{ p}\prod_{ l=1}^{ p}\|
f_{l}\|_{ pq/( ( p-1)q+1)}\label{60.1}
\end{equation}
and for any $0<\si<d$
\begin{equation}
\Bigg|\int_{ (R^{ d})^{ p}} {\prod_{ l=1}^{ p}f_{l}( x_{ l})\over |x_{ 1}+\cdots+x_{
p}|^{ d-\si}}
\prod_{ j=1}^{ p}\,d x_{ l}\Bigg|\leq  C^{ p}\prod_{ l=1}^{ p}\|
f_{l}\|_{ pd/( ( p-1)d+\si)}.\label{60.2}
\end{equation}

Furthermore, for any $n$ and any $F_{l}=F_{l}(x_{l, j};\,1\leq j\leq n)$, $1\leq l\leq p$ we have 
\bea&&
\Bigg|\int_{ (R^{ d})^{ np}} \prod_{j=1}^{ n}{1\over |x_{1,j}+\cdots+x_{
p,j}|^{ d-\si}}\,\,\prod_{ l=1}^{ p}F_{l}\,\,\prod_{j=1}^{ n} \prod_{l=1}^{ p}\,d x_{
l,j}\Bigg|\label{60.2n}\\
&&\leq  C^{ p}\prod_{ l=1}^{ p}
\|F_{l}\|_{ pd/( ( p-1)d+\si)}.\nn 
\eea
and more generally,  for some $C<\ff$ independent of $z\in R^{d}$
\bea&&
\Bigg|\int_{ (R^{ d})^{ np}} \prod_{j=1}^{ n}{1\over |x_{1,j}+\cdots+x_{
p,j}-z|^{ d-\si}}\,\,\prod_{ l=1}^{ p}F_{l}\,\,\prod_{j=1}^{ n} \prod_{l=1}^{ p}\,d x_{
l,j}\Bigg|\label{60.2nz}\\
&&\leq  C^{ p}\prod_{ l=1}^{ p}
\|F_{l}\|_{ pd/( ( p-1)d+\si)}.\nn 
\eea
\end{lemma}

\Proof  We prove (\ref{60.1}) by induction on $p$. The case $p=1$ is trivial. Thus
assume (\ref{60.1}) holds for all $p\leq m-1$. Since $t^{ -1}=r^{ -1}+s^{ -1}-1$
when $t=q, r=mq/( ( m-1)q+1),s=mq/(  m-1+q)$, it follows from Young's inequality, 
 \cite{D}, p. 275, that
\begin{equation}
\|f_{ 1}\ast\cdots\ast f_{m}\|_{ q}\leq C\|f_{1}\|_{ mq/( ( m-1)q+1)}
\|f_{ 2}\ast\cdots\ast f_{m}\|_{mq/(  m-1+q)}.\label{60.3}
\end{equation}
By our induction hypothesis and using the fact that
\begin{equation}
{(m-1)mq \over (m-2)mq+m-1+q}={(m-1)mq \over (m-1)^{2}q+m-1}
={mq \over (m-1)q+1}\label{frac.1}
\end{equation}
we see that
\begin{equation}
\|f_{ 2}\ast\cdots\ast f_{m}\|_{mq/(  m-1+q)}\leq C^{ m-1}\prod_{l=2}^{ m}\|
f_{l}\|_{ mq/( ( m -1)q+1)}\label{60.4}
\end{equation}
which completes the proof of (\ref{60.1}).

To prove (\ref{60.2}) we write
\begin{equation}
\int_{ (R^{ d})^{p}} {\prod_{ l=1}^{ p}f_{l}( x_{ l})\over |x_{ 1}+\cdots+x_{
p}|^{ d-\si}}
\prod_{l=1}^{ p}\,d x_{ l}=
\int_{ (R^{ d})^{2}} {f_{1}( x)\,\,\, (f_{2}\ast\cdot\ast f_{p})( y)\over |x-y|^{ d-\si}}
\,d x\,dy\label{frac.2}
\end{equation}
and apply (\ref{sobolineq}) with $r=pd/( ( p-1)d+\si)$ so that $s=pd/( ( p-1)\si+d)$ to obtain
\begin{eqnarray} &&
\Bigg|\int_{ (R^{ d})^{ p}} {\prod_{ l=1}^{ p}f_{l}( x_{ l})\over |x_{ 1}+\cdots+x_{
p}|^{ d-\si}}
\prod_{l=1}^{ p}\,d x_{l}\Bigg|\label{60.2m}\\
&&\leq  C \|f_{1}\|_{ pd/( ( p-1)d+\si)}\,\,\,
\|f_{2}\ast\cdot\ast f_{p}\|_{pd/( ( p-1)\si+d)}.\nn
\end{eqnarray}
Then using  (\ref{60.1}) and  the fact that
\begin{equation}
{(p-1)pd \over (p-2)pd+(p-1)\si+d}={(p-1)pd\over (p-1)^{2}d+(p-1)\si}
={pd \over (p-1)d+\si}\label{frac.3}
\end{equation}
we obtain (\ref{60.2}).

We next prove (\ref{60.2n}). By (\ref{60.2})
\bea&&\qquad
\Bigg|\int_{ (R^{ d})^{ np}} \prod_{j=1}^{ n}{1\over |x_{1,j}+\cdots+x_{
p,j}|^{ d-\si}}\,\,\prod_{ l=1}^{ p}F_{l}\,\prod_{j=1}^{ n} \prod_{l=1}^{ p}\,d x_{
l,j}\Bigg|\label{60.6}\\
&&\leq \int_{ (R^{ d})^{ (n-1)p}} \prod_{j=2}^{ n}{1\over |x_{1,j}+\cdots+x_{
p,j}|^{ d-\si}}\nn\\
&&\qquad\Bigg|\int_{ (R^{ d})^{ p}} {\prod_{ l=1}^{ p}|F_{l}|\over |x_{
1,1}+\cdots+x_{ p,1}|^{ d-\si}}
\prod_{l=1}^{ p}\,d x_{l,1}\Bigg|\,\,\,\prod_{j=2}^{ n} \prod_{l=1}^{ p}\,d x_{
l,j}\nn \\
&&\leq \int_{ (R^{ d})^{  (n-1)p}} \prod_{j=2}^{ n}{1\over |x_{1,j}+\cdots+x_{
p,j}|^{ d-\si}}\nn \\
&&\hspace{ 1in}
\prod_{ l=1}^{ p}
\|F_{l}\|_{pd/( ( p-1)d+\si), x_{l,1}}
\prod_{j=2}^{ n} \prod_{l=1}^{ p}\,d x_{
l,j}\nn
\eea
where
 \[\|F_{l}\|_{ q,\,x_{l,1}}=\(\int_{ R^{ dp}}
|F_{l}(x_{l, j};\,1\leq j\leq n)|^{ q}\,dx_{l,1 }\)^{ 1/q}.\]
 Inequality (\ref{60.2n}) then
follows by iterating this step. For example, the next iteration will bound (\ref{60.6}) 
by 
\bea&&
 \int_{ (R^{ d})^{  (n-2)p}} \prod_{j=3}^{ n}{1\over |x_{1,j}+\cdots+x_{
p,j}|^{ d-\si}}\label{60.7}\\
&&\hspace{ 1in} 
\prod_{ l=1}^{ p}
\|F_{l}\|_{pd/( ( p-1)d+\si), x_{l,1},  x_{l,2}}
\prod_{j=3}^{ n} \prod_{l=1}^{ p}\,d x_{
l,j}\nn\nn
\eea
where now 
\begin{eqnarray} &&
\|F_{l}\|_{q, x_{l,1},  x_{l,2}}\label{60.8}\\ &&
=\(\int_{
R^{ d}}  \|F_{l}\|_{ q,\,x_{l,1}}^{ q}\,dx_{l,2}\)^{ 1/q}\nonumber\\
 &&
=\(\int_{
R^{2 d}}|F_{l}(x_{l, j};\,1\leq j\leq n)|^{ q}\,dx_{l,1 }\,dx_{l,2 }\)^{ 1/q}.\nonumber
\end{eqnarray}
 It should 
be clear that this will lead to (\ref{60.2n}).

Now let $T^{l,j}_{z}$ denote translation of $x_{l,j}$ by $z$ and set 
$\mathcal{T}=\prod_{j=1}^{n}T^{1,j}_{z}$. Then using (\ref{60.2n}) and the translation invariance of Lebesgue measure
\bea&&
\Bigg|\int_{ (R^{ d})^{ np}} \prod_{j=1}^{ n}{1\over |x_{1,j}+\cdots+x_{
p,j}-z|^{ d-\si}}\,\,\prod_{ l=1}^{ p}F_{l}\,\,\prod_{j=1}^{ n} \prod_{l=1}^{ p}\,d x_{
l,j}\Bigg|\label{60.2nz}\\
&&=\Bigg|\int_{ (R^{ d})^{ np}} \prod_{j=1}^{ n}{1\over |x_{1,j}+\cdots+x_{
p,j}|^{ d-\si}}\,\,\prod_{ l=1}^{ p}\mathcal{T}F_{l}\,\,\prod_{j=1}^{ n} \prod_{l=1}^{ p}\,d x_{
l,j}\Bigg|\nn\\
&&\leq  C^{ p} \prod_{l=1}^{ p}\|
\mathcal{T}F_{l}\|_{ pd/( ( p-1)d+\si)}=C^{ p} \prod_{l=1}^{ p}\|
F_{l}\|_{ pd/( ( p-1)d+\si)}\nn 
\eea
which is (\ref{60.2nz}).
\qed

\def\noopsort#1{} \def\printfirst#1#2{#1} \def\singleletter#1{#1}
  \def\switchargs#1#2{#2#1} \def\bibsameauth{\leavevmode\vrule height .1ex
  depth 0pt width 2.3em\relax\,}
\makeatletter \renewcommand{\@biblabel}[1]{\hfill#1.}\makeatother
\newcommand{\bysame}{\leavevmode\hbox to3em{\hrulefill}\,}

\end{document}